\documentclass[twoside]{amsart}

\usepackage{array}
\usepackage[dvips]{color}
\usepackage{amsfonts}
\usepackage{amssymb, latexsym, amsmath, pb-diagram}
\usepackage{pstricks-add}
\usepackage{lamsarrow, pb-lams}
\usepackage{multicol}
\usepackage{bm}
\usepackage[mathscr]{eucal}
\usepackage[all]{xy}
\usepackage{epic,eepic}

\allowdisplaybreaks

\def\q{\hfill q.e.d.}
\def\qb{\hfill $\Box$}

\def\cf{{\it cf.}\ }

\begin{document}

\title{On the rational cohomology  of spin hyperelliptic mapping class groups}

\author{Gefei Wang}

\email{gefeiw@ms.u-tokyo.ac.jp}

\address{Graduate School of Mathematical Science, University of Tokyo, 3-8-1 Komaba, Meguro-ku, Tokyo, 153-8914, Japan}

\thanks{This project is supported by NSFC No.11871284}



\begin{abstract}
Let $\mathfrak{G}$ be the subgroup $\mathfrak{S}_{n-q} \times \mathfrak{S}_{q}$ of the $n$-th symmetric group $\mathfrak{S}_{n}$ for $n-q \geq q$. In this paper, we study the  $\mathfrak{G}$-invariant part of the rational cohomology group of the pure braid group $P_{n}$. The invariant part includes the rational cohomology of a spin hyperelliptic mapping class group of genus $g$ as a subalgebra when $n=2g+2$, denoted by $H^*(P_{n})^{\mathfrak{G}}$. Based on the study of Lehrer-Solomon \cite{LS86},
we prove that they are independent of $n$ and $q$ in degree $*\leq q-1$. We also give a formula to calculate the dimension of $H^*(P_{n})^{\mathfrak{G}}$ and calculate it in all degree for $q\leq 3$.
\end{abstract}

\maketitle

\section{Introduction}
Let $\Sigma_g$ be a compact connected oriented surface of genus $g$, which admits a branched $2$-fold covering $\varphi: \Sigma_g \rightarrow S^2$ with $2g+2$ branch points.
The branched covering induces a group homomorphism of the Artin braid group $B_{2g+2}$ to the mapping class group of genus $g$.

A $2$-Spin structure $c$, or spin structure $c$ for short, on the surface $\Sigma_g$
can be regarded as a map $$c :  H_1(\Sigma_g , \mathbb{Z}/2) \rightarrow \mathbb{Z}/2$$ such that $c(x+y)=c(x)+c(y)+x \cdot y$ for any $x,y \in H_1(\Sigma_g , \mathbb{Z}/2)$, where $x\cdot y$ is the intersection number. The mapping class group of genus $g$ acts on the set of spin structures in a natural way, and we call the isotropy group at a spin structure $c$ the spin mapping class group associated with $c$ ({\it cf.} \cite{J80}).

By the classical theory of Riemann surfaces, the set of spin structures
is naturally bijective to the set of subsets of the branch points of the covering $\varphi$ whose cardinality is congruent to $g+1$ modulo $2$ ({\it cf.} e.g., \cite{M84} Proposition 3.1, p.3.95).  The bijection is equivariant under the action
of $B_{2g+2}$. Hence the action of the Artin braid group is reduced to that of the symmetric group $\mathfrak{S}_{2g+2}$, and the isotropy group of each spin structure is isomorphic to that of the corresponding subset of the set of branch points. The action has $\lceil \frac{g}{2}\rceil+1$ orbits.
 The isotropy group $\mathfrak{G}_i$ of each orbits can be described by the direct products of two symmetric groups $\mathfrak{S}_{g+1+2i} \times \mathfrak{S}_{g+1-2i}$ or a $\mathbb{Z}/2$-extension of direct products of two symmetric groups $\mathfrak{S}_{g+1} \overleftrightarrow{\times} \mathfrak{S}_{g+1}$. The author also gives an alternative proof of these results in a purely combinatorial way in \cite{W23}.

In \cite{RM14}, Randal-Williams established the homology stability of the moduli space of $r$-Spin surfaces with underlying surface $\Sigma_{g,b}$ and fixed $r$-Spin structure $c$ along the boundary, which is called the stable $r$-Spin mapping class group, denoted by $\mathcal{M}^{Spin^r}(\Sigma_{g,b};c)$, where $\Sigma_{g,b}$ is the compact connected oriented surface of genus $g$ with $b$ boundary components. In particular, for $r=2$, the case of ordinary spin Riemann surfaces, the integral homology groups $$H_*(\mathcal{M}^{Spin^2}(\Sigma_{g,b};c))$$ are independent of $g$, $b$ and $\delta$ in degree $5*\leq 2g-7$.
Randal-Williams also computed the rational cohomology of $\mathcal{M}^{Spin^r}(\Sigma_{g,b};c)$. In Theorem $1.3$ \cite{RM12}, he proved the rational cohomology algebra is freely generated by the Mumford-Morita-Miller classes.

Now we want to study unstable aspects of the rational cohomology of the spin mapping class group. As first examples, we have spin hyperelliptic mapping class groups.
Let $\mathcal{S}(\Sigma_g)$ be the hyperelliptic mapping class group for a compact connected oriented surface $\Sigma_g$ of genus $g$.
 We can consider the spin hyperelliptic mapping class group $\mathcal{S}(\Sigma_g ; c)$ as the isotropy group of $\mathcal{S}(\Sigma_g)$ at each spin structure c. As was stated above,  the action of the group $\mathcal{S}(\Sigma_g)$ on the set of spin structures on the surface $\Sigma_g$ has $\lceil \frac{g}{2}\rceil+1$ orbits. Hence we have  $\lceil \frac{g}{2}\rceil+1$ spin hyperelliptic mapping class groups for $\Sigma_g$. Now we consider the rational cohomology groups $H^*(\mathcal{S}(\Sigma_{g};c))$.
  As is known classically, the Mumford-Morita-Miller classes vanish on the rational cohomology of $\mathcal{S}(\Sigma_g)$.

   In \cite{BH71}, Birman and Hilden proved that the mapping class group of the $2$-sphere with $2g+2$ branch points, denoted by $\mathcal{M}(\Sigma_{0,2g+2})$, is isomorphic to the hyperelliptic mapping class group of $\Sigma_g$ modulo a subgroup generated by an involution $\iota$ of order $2$. The mapping class group $\mathcal{M}(\Sigma_{0,2g+2})$ has a natural surjective homomorphism onto the $(2g+2)$-th symmetric group $\mathfrak{S}_{2g+2}$. Hence it induces surjective homomorphisms of $\mathcal{S}(\Sigma_{g})$ and the $(2g+2)$-th sphere braid group onto $\mathfrak{S}_{2g+2}$. We denote the kernels
	of these three homomorphisms by $\mathcal{M}^P(\Sigma_{0,2g+2})$, $\mathcal{S}^P(\Sigma_{g})$ and $SP_{2g+2}$, respectively. The kernel $SP_{2g+2}$ coincides with the $(2g+2)$-th pure sphere braid group. The group $\mathcal{S}^P(\Sigma_g)$ is the quotient of $SP_{2g+2}$ by a normal subgroup of order $2$. In general, if a group $Q$ is the quotient of a group $G$ by a normal subgroup of order $2$, then the rational cohomology of $Q$ is isomorphic to that of $G$. Hence the rational cohomology groups of $\mathcal{M}^P(\Sigma_{0,2g+2})$, $\mathcal{S}^P(\Sigma_{g})$ and $SP_{2g+2}$ are $\mathfrak{S}_{2g+2}$-equivariantly isomorphic to each other. Thus  $H^*(\mathcal{S}(\Sigma_{g};c))$ is isomorphic to the  $\mathfrak{G}_i$-invariant part of the rational cohomology group of the pure sphere braid group $SP_{2g+2}$, denoted by $H^*(SP_{2g+2})^{\mathfrak{G}_i}$. Moreover, for any $n \geq 2$, the natural surjection of the pure Artin braid group $P_n$ onto the pure sphere braid group $SP_n$ induces an injection of the rational cohomology groups (\cf e.g., \cite{K90}). Hence
    the  $\mathfrak{G}_i$-invariant part of the rational cohomology group of  $P_{2g+2}$, denoted by $H^*(P_{2g+2})^{\mathfrak{G}_i}$, includes $H^*(\mathcal{S}(\Sigma_{g};c))$ as a subalgebra.

 In this paper, we study  $H^*(P_{n})^{\mathfrak{G}}$, where $\mathfrak{G}=\mathfrak{S}_{n-q} \times \mathfrak{S}_{q}$, $n-q \geq q$, as a first step.
$H^*(P_{n})^{\mathfrak{G}}$ equals the rational cohomology group of the isotropy subgroup in the Artin braid group $B_{n}$ for the orbit whose isotropy group is $\mathfrak{G}$.

In section $2$,  we give a description of $H^*(P_{n})^{\mathfrak{G}}$ based on \cite{LS86} in Lemma $2.1$. More explicitly, as a $\mathbb{Q}$-vector space, we have $$H^*(P_{n})^{\mathfrak{G}}\cong\bigoplus_{\lambda}\bigoplus_{(s \otimes \zeta_{\lambda},1)_{(Z_{\lambda})_s}=1 \atop s\in \mathfrak{G}\backslash\mathfrak{S}_n / Z_{\lambda}}(s \otimes \zeta_{\lambda}),$$
 where $\lambda$ are the $(n-*)$ partitions of $n$, $\zeta_{\lambda}$ is a $1$-dimensional representation of $\lambda$, $Z_{\lambda}$ is the centralizer of $\lambda$,  $(Z_{\lambda})_s=sZ_{\lambda}s^{-1} \cap \mathfrak{G}$  and $s \otimes \zeta_{\lambda}$ is defined by $(s \otimes \zeta_{\lambda})(g)=\zeta_{\lambda}(s^{-1}gs)$.
  We give a complete invariant $$\chi: \mathfrak{G}\backslash\mathfrak{S}_n / Z_{\lambda}\rightarrow \{\emptyset\}\sqcup(\coprod_{d}(\mathbb{Z}_{\geq0})^{d}/\langle c_{d}\rangle)$$ for $H^*(P_{n})$ in Definition $2.4$, where $c_d$ is a cycle of length $d$. In Theorem $2.5$, we give an equivalence between the pair $(\lambda,s)$ with $(s \otimes \zeta_{\lambda},1)_{(Z_{\lambda})_s}=1$ and the invariant $\chi$.
We calculate a ring structure of $H^*(P_{n})^{\mathfrak{G}}$ in Theorem $2.6$.
Let $\mathfrak{G}_1=\mathfrak{S}_{n-q}\times\mathfrak{S}_{q}$ and $\mathfrak{G}_2=\mathfrak{S}_{n-q-1}\times\mathfrak{S}_{q+1}$, $n-q-1\geq q+1$, we prove that there exists a bijection
\[  \psi: H^{*}(P_n)^{\mathfrak{G}_{1}}
\rightarrow H^{*}(P_n)^{\mathfrak{G}_{2}}\]
in the range $*\leq q-1$ in Theorem $2.7$.

 Section $3$ contains a formula to calculate the dimension of $H^*(P_{n})^{\mathfrak{G}}$.

{\noindent\bf Theorem 3.1.}
{\it The dimension of $H^{*}(P_n)^\mathfrak{G}$ equals
$$
\sum_{(\lambda,d)\in M}\prod_{i\in M_1(\lambda,d)}
(\sum_{b=1}^{|(\lambda_i,d_i)|}|K( b)|\left(\begin{array}{c}|\Pi(\lambda_i,d_i)|\\b\end{array}\right))
\prod_{i\in M_2(\lambda,d)}\left(\begin{array}{c}|\Pi(\lambda_i,d_i)|\\|(\lambda_i,d_i)|\end{array}\right),
$$
where
$\left(\begin{array}{c}|\Pi(\lambda_i,d_i)|\\b\end{array}\right)$, $1\leq b\leq |\Pi(\lambda_i,d_i)|$, is the binomial coefficient and
 we write
$\left(\begin{array}{c}|\Pi(\lambda_i,d_i)|\\b\end{array}\right)=0$
if $|\Pi(\lambda_i,d_i)|< b$.}

\noindent Here we denote $\Lambda(n-*)$  the set of all $(n-*)$ partitions of $n$. For each partition $\lambda=(\lambda_1,\lambda_2,\cdots,\lambda_a, 1,\cdots,1)\in\Lambda(n-*)$, $\lambda_i\geq 2$, $1\leq i \leq a$,
 we define $$d=(d_1,d_2,\cdots,d_a)\in(\mathbb{Z}_{\geq0})^a$$ and consider
$$(\lambda,d)=((\lambda_1,d_1),(\lambda_2,d_2),\cdots,(\lambda_a,d_a))$$ who satisfies $d_i\geq d_{i+1}$ if $\lambda_i=\lambda_{i+1}$.
We denote
 \begin{align*}
M&=\{(\lambda,d)|\lambda=(\lambda_1,\lambda_2,\cdots,\lambda_a,1,1,\cdots)\in\Lambda(n-*),\sum_{i=1}^{a}d_i\leq q\}\\
M_1(\lambda,d)&=\{1\leq i \leq a|  \lambda_i\equiv 1 \text{ mod } 2\},\\
M_2(\lambda,d)&=\{1\leq i \leq a|\lambda_i \equiv 0 \text{ mod } 2\},\\
|(\lambda_i,d_i)|&=\sharp\{1\leq j\leq a|(\lambda_j,d_j)=(\lambda_{i},d_i)\}
\end{align*}
and
\begin{align*}
K( b)&=\{(k_1,k_2,\cdots,k_b)\in (\mathbb{Z}_{>0})^b| k_1+k_2+\cdots+k_b=|(\lambda_i,d_i)|\}
\end{align*}
for $d_i\in \mathbb{Z}_{\geq0}$, $1\leq b \leq |(\lambda_i,d_i)|$. The cardinality $|\Pi(\lambda_i,d_i)|$ will be given in Theorem $3.3$.
Furthermore, Theorem $3.3$ shows that $|\Pi(\lambda_i,d_i)|$ is only depending on $\lambda_i$ and $d_i$. From Theorem $2.7$, Theorem $3.1$ and Theorem $3.3$, since the proof of Theorem $2.7$ shows that $d_1+d_2+\cdots+d_a\leq *+1$, we have $H^*(P_{n})^{\mathfrak{G}}$ are independent of $n$ and $q$ in degree $*\leq q-1$.
We also calculate the dimension of $H^*(P_{n})^{\mathfrak{G}}$ in all degree for $q\leq 3$.
More explicitly, in Theorem $3.4$, we prove that
$$H^{i}(P_n)^{\mathfrak{G}}=\left\{\begin{array}{ccc}\mathbb{Q}, & \text{if }i=0,n-1,\\
\mathbb{Q}\oplus\mathbb{Q}, & \text{if }1\leq i \leq n-2, \end{array}\right.$$
when $q=1$.
In Theorem $3.5$,
$$H^{i}(P_n)^{\mathfrak{G}}=\left\{\begin{array}{ccc}
\mathbb{Q}, & \text{if }i=0, &\\
\mathbb{Q}^{2i}, & \text{if }i \equiv 0,2,& 1\leq i \leq n-3 , \\
\mathbb{Q}^{2i+1}, &   \text{if }i\equiv 1,&  1\leq i\leq n-3, \\
\mathbb{Q}^{2i-1}, &  \text{if }i\equiv 3,& 1\leq i \leq n-3, \\
\mathbb{Q}^{\frac{3}{2}i}, & \text{if }i \equiv 0,2, & i=n-2,\\
\mathbb{Q}^{\frac{3}{2}i-\frac{1}{2}}, & \text{if }i \equiv 1, & i=n-2,\\
\mathbb{Q}^{\frac{3}{2}i+\frac{1}{2}}, & \text{if }i \equiv 3, & i=n-2,\\
\mathbb{Q}^{\frac{1}{2}i},  & \text{if }i \equiv 0,2, & i=n-1,\\
\mathbb{Q}^{\frac{1}{2}i+\frac{1}{2}},  & \text{if }i \equiv 1, & i=n-1,\\
\mathbb{Q}^{\frac{1}{2}i-\frac{1}{2}},  & \text{if }i \equiv 3, & i=n-1,
\end{array}\right.\text{ mod }4,$$
 when $q=2$.
Finally in  Theorem $3.6$,
$$H^{i}(P_n)^{\mathfrak{G}}=\left\{\begin{array}{ccc}
\mathbb{Q}^{i^2-\frac{1}{3}i+1}, & \text{if }i \equiv 0,3, & 0\leq i\leq n-4,\\
\mathbb{Q}^{i^2-\frac{1}{3}i+\frac{7}{3}}, & \text{if }i \equiv 1,10, & 0\leq i\leq n-4,\\
\mathbb{Q}^{i^2-\frac{1}{3}i+\frac{5}{3}}, & \text{if }i \equiv 2,5, & 0\leq i\leq n-4,\\
\mathbb{Q}^{i^2-\frac{1}{3}i+\frac{4}{3}},  & \text{if }i \equiv 4,7, & 0\leq i\leq n-4,\\
\mathbb{Q}^{i^2-\frac{1}{3}i+2,}  & \text{if }i \equiv 6,9, &0\leq i\leq n-4\\
\mathbb{Q}^{i^2-\frac{1}{3}i+\frac{2}{3}},  & \text{if }i \equiv 8,11, & 0\leq i\leq n-4,\\
\mathbb{Q}^{\frac{11}{12}i^2-\frac{5}{12}i+1},  & \text{if }i\equiv 0,3,4,7, &i=n-3,\\
\mathbb{Q}^{\frac{11}{12}i^2-\frac{5}{12}i+\frac{3}{2}},  & \text{if }i\equiv 1,6,9,10, &i=n-3,\\
\mathbb{Q}^{\frac{11}{12}i^2-\frac{5}{12}i+\frac{7}{6}},  & \text{if }i\equiv 2,5,&i=n-3,\\
\mathbb{Q}^{\frac{11}{12}i^2-\frac{5}{12}i+\frac{2}{3}},  & \text{if }i\equiv 8,11,&i=n-3,\\
\mathbb{Q}^{\frac{7}{12}i^2-\frac{5}{12}i},  & \text{if }i\equiv 0,3,8,11,&i=n-2,\\
\mathbb{Q}^{\frac{7}{12}i^2-\frac{5}{12}i+\frac{5}{6}},  & \text{if }i\equiv 1,10,&i=n-2,\\
\mathbb{Q}^{\frac{7}{12}i^2-\frac{5}{12}i+\frac{1}{2}},  & \text{if }i\equiv 2,5,6,9,&i=n-2,\\
\mathbb{Q}^{\frac{7}{12}i^2-\frac{5}{12}i+\frac{1}{3}},  & \text{if }i\equiv 4,7,&i=n-2,\\
\mathbb{Q}^{\frac{1}{6}i^2-\frac{1}{6}i},  & \text{if }i \equiv 0,1,3,4,6,7,9,10,& i=n-1,\\
\mathbb{Q}^{\frac{1}{6}i^2-\frac{1}{6}i-\frac{1}{3}},  & \text{if }i \equiv 2,5,8,11,& i=n-1,\\
\end{array}\right.\text{ mod }12,$$
when $q=3$.

The author is grateful to Professor Nariya Kawazumi for his encouragement during this
research.

\section{$\mathfrak{G}_i$-invariant part of $H^*(P_{n})$}
In this section, we compute the $\mathfrak{G}$-invariant part of the rational cohomology $H^*(P_{n})$.

 By definition,  the pure braid group $P_n$ is the fundamental group of the configuration space
\[X_n := \{(u_1, \cdots, u_n) \in \mathbb{C}^n\mid u_i \neq u_j \text{ for } i\neq j\}.\]
From the study by Arnold \cite{A69}, the cohomology ring $H^*(P_{n})$ is the quotient of an exterior graded ring, denoted by $A(n)$, generated by
$\left(
                                 \begin{array}{c}
                                   n \\
                                   2 \\
                                 \end{array}
                               \right)$ differential forms $\omega_{i,j}=\omega_{j,i}=\frac{du_i-du_j}{u_i-u_j}$, $u_i \in \mathbb{C}$, $1\leq i \neq j \leq n$, modulo the $\left(
                                 \begin{array}{c}
                                   n \\
                                   3 \\
                                 \end{array}
                               \right)$ relations
                               \[\omega_{k.l}\omega_{l,m}+\omega_{l,m}\omega_{m,k}+\omega_{m,k}\omega_{k,l}=0,\]
where $1\leq k\neq l \neq m\leq n$.

The symmetric group $\mathfrak{S}_n$ acts on $A(n)$ in a natural way.
 For a partition $\lambda=(\lambda_{1}\geq\lambda_{2}\geq\cdots\geq\lambda_{j}>0)$ of $n$, let $c_{\lambda}=g_1g_2\cdots g_j$ be a element of $\mathfrak{S}_n$, where
 $$g_i=(\lambda_1+\lambda_2+\cdots+\lambda_{i-1}+1,\lambda_1+\lambda_2+\cdots+\lambda_{i-1}+2,\cdots,\lambda_1+\lambda_2+\cdots+\lambda_{i})$$ is a cycle of length $\lambda_i$.
Let $N_{\lambda}$ be the subgroup of $\mathfrak{S}_n$ generated by $\{\nu_i|\lambda_i=\lambda_{i+1}\}$, where
 \begin{align*}
 \nu_i=&(\lambda_1+\lambda_2+\cdots+\lambda_{i-1}+1,\lambda_1+\lambda_2+\cdots+\lambda_{i}+1)\\
 &(\lambda_1+\lambda_2+\cdots+\lambda_{i-1}+2,\lambda_1+\lambda_2+\cdots+\lambda_{i}+2)\cdots\\
 &(\lambda_1+\lambda_2+\cdots+\lambda_{i},\lambda_1+\lambda_2+\cdots+\lambda_{i+1}),
 \end{align*}
 $C_{\lambda}=\langle g_1\rangle\times \langle g_2\rangle\times\cdots\times \langle g_p\rangle$  the direct product of cyclic groups $\langle g_i\rangle$.
 Let $Z_{\lambda}$ be the centralizer of the element $c_{\lambda}$.
 Then the  centralizer $Z_{\lambda}$ of the element $c_{\lambda}$ is the semi-product of $N_{\lambda}$ and $C_{\lambda}$, $Z_{\lambda}=N_{\lambda}\ltimes C_{\lambda}$.

  Assume that $\alpha_{\lambda}$ is the character of $N_{\lambda}$ defined by $\alpha_{\lambda}(\nu_i)=(-1)^{\lambda_i+1}$, $\phi_{\lambda}=(\phi_{\lambda_1}\otimes\phi_{\lambda_2}\otimes\cdots\otimes\phi_{\lambda_{p}})\varepsilon$ is the character of $C_{\lambda}$, where $\phi_{\lambda_i}(g_i)=e^{\frac{2\pi i}{\lambda_i}}$ and $\varepsilon$ is the sign character of $\mathfrak{S}_n$.
 From \cite{LS86}, the $i$-th cohomology group of $P_n$ has a $\mathbb{Q}\mathfrak{S}_n$-module isomorphism
\[ H^i(P_{n}) \cong \bigoplus_{\lambda}\text{Ind}_{Z_{\lambda}}^{\mathfrak{S}_n}(\zeta_{\lambda}),\]
where $\lambda=(\lambda_{1}\geq\lambda_{2}\geq\cdots\geq\lambda_{j}>0)$ are those $j$ partitions of $n$ such that $i+j=n$ and $\zeta_{\lambda}=\alpha_{\lambda}\phi_{\lambda}$.

Let $\rho$ be the representation of $\mathfrak{S}_n$ on $H^*(P_{n})$.
The dimension of $\mathfrak{S}_n$-invariant part of $H^*(P_{n})$ can be computed by $$(\rho,1)=\frac{1}{|\mathfrak{S}_n|}\sum_{\sigma\in\mathfrak{S}_n}\rho(\sigma)\overline{1}(\sigma).$$ By using Frobenius reciprocity,
\begin{align*}
(\rho,1)=&\sum_{\lambda}(1,\zeta_{\lambda})_{Z_{\lambda}}\\
=&\sum_{\lambda}\frac{1}{|Z_{\lambda}|}\sum_{z\in Z_{\lambda}}1(z)\overline{\zeta_{\lambda}}(z)\\
=&2
\end{align*}
since $(1,\zeta_{\lambda})_{Z_{\lambda}}=1$ if and only if $\lambda=(1,1,\cdots,1)$ or $\lambda=(2,1,1,\cdots,1)$. Namely, only the $0$-th and the first cohomology groups have none-trivial $\mathfrak{S}_n$-invariant part (cf. \cite{A69}).
  The  $\mathfrak{S}_n$-invariant part in degree $1$ is generated by the total sum of the generators $\omega_{i,j}$ of $H^1(P_{n})$ and those of the other degrees are trivial.

Let $\mathfrak{G}$ be a subgroup of $\mathfrak{S}_{n}$.
Since, as a $\mathbb{Q}\mathfrak{G}$-module, \[ H^*(P_{n}) \cong \bigoplus_{\lambda}\text{Res}_{\mathfrak{G}}\text{Ind}_{Z_{\lambda}}^{\mathfrak{S}_n}(\zeta_{\lambda})
\cong\bigoplus_{\lambda}\bigoplus_{s\in \mathfrak{G}\backslash\mathfrak{S}_n / Z_{\lambda}}\text{Ind}_{(Z_{\lambda})_s}^{\mathfrak{G}}(s \otimes \zeta_{\lambda}),\]
where $(Z_{\lambda})_s=sZ_{\lambda}s^{-1} \cap \mathfrak{G}$ and $s \otimes \zeta_{\lambda}$ is defined by $(s \otimes \zeta_{\lambda})(g)=\zeta_{\lambda}(s^{-1}gs)$,
hence one can compute the $\mathfrak{G}$-invariant part of $H^*(P_{n})$ by $(s \otimes \zeta_{\lambda},1)_{(Z_{\lambda})_s}$.

{\noindent\bf Lemma 2.1.}
{\it For a partition $\lambda=(\lambda_{1}\geq\lambda_{2}\geq\cdots\geq\lambda_{j}>0)$ of $n$ and an element $s$
   of $\mathfrak{G}\backslash\mathfrak{S}_n \slash Z_{\lambda}$, we have

 $$(s \otimes \zeta_{\lambda},1)_{(Z_{\lambda})_s}=
 \left\{\begin{array}{cc} 0, & \exists z=\nu c\in Z_{\lambda}$\text{ s.t. }$  szs^{-1} \in \mathfrak{G}, \zeta_{\lambda}(z)=\alpha_{\lambda}(\nu)\phi_{\lambda}(c)\neq1, \\
 1, & \text{otherwise.} \end{array}\right.$$}

 {\bf Proof:} Since $s \otimes \zeta_{\lambda}$ is a 1-dimensional representation,$(s \otimes \zeta_{\lambda},1)_{(Z_{\lambda})_s}$ equals 0 or 1, and the
equation $s \otimes \zeta_{\lambda}=0$ is equivalent to the representation $s \otimes \zeta_{\lambda}$ is non-trivial. This
proves the lemma.\qb

From Lemma $2.1$, as a $\mathbb{Q}$-vector space, we have $$H^*(P_{n},\mathbb{Q})^{\mathfrak{G}}\cong\bigoplus_{\lambda}\bigoplus_{(s \otimes \zeta_{\lambda},1)_{(Z_{\lambda})_s}=1 \atop s\in \mathfrak{G}\backslash\mathfrak{S}_n / Z_{\lambda}}(s \otimes \zeta_{\lambda}).$$

Firstly, we consider the case $\mathfrak{G}=\mathfrak{S}_{n-q} \times \mathfrak{S}_q$, $n-q\geq q$. Denote by $[n]$ the index set $$\{1,2,\cdots, n\},$$ by
$P(i)=\{A\subset [n]\}$ the set of subsets of $[n]$, where $|A|=i$, by $\mathfrak{P}$ the index set $$\{1,2,\cdots,n-q\}$$ and by $\mathfrak{Q}$ the index set $$\{n-q+1,n-q+2,\cdots,n\}.$$
 Then we have $[n]=\mathfrak{P} \cup \mathfrak{Q}$.

{\noindent\bf Lemma 2.2.}
{\it For $\mathfrak{G}=\mathfrak{S}_{n-q} \times \mathfrak{S}_q$, $n-q\geq q\geq 1$, we can take the set
$$S:=\{s=\prod^{d}_{i=1}(k_i,l_i)|0\leq d\leq q, k_i<k_{i+1}\in\mathfrak{P},l_i<l_{i+1}\in\mathfrak{Q} \}$$ as a set of representatives of $\mathfrak{G}\backslash \mathfrak{S}_{n}$.}

{\noindent\bf Proof.}  Let $$P(n-q,q)=\{(A,B)\in P(n-q)\times P(q)|A\cap B=\emptyset\}\subset P(n-q)\times P(q).$$
  $\mathfrak{S}_{n}$ acts on $P(n-q,q)$ transitively. Since
 $\mathfrak{G}$ is the isotopy group of $\mathfrak{P}\times\mathfrak{Q}$, we have $|\mathfrak{G}\backslash \mathfrak{S}_{n}|=|P(n-q,q)|$. On the other hand, we have
 $S \circ (\mathfrak{P}\times\mathfrak{Q})=P(n-q,q)$.
We compute the cardinality of $S$ and $\mathfrak{G}\backslash \mathfrak{S}_{n}$.
Since
\begin{align*}
|S|=& \sum_{d=0}^{q}\left(\begin{array}{c}n-q\\d\end{array}\right)\left(\begin{array}{c}q\\d\end{array}\right)\\
=&\left(\begin{array}{c}n\\q\end{array}\right)\\
=&|\mathfrak{G}\backslash \mathfrak{S}_{n}|,
\end{align*}
the lemma follows.\qb

{\noindent\bf Lemma 2.3.}
{\it  Let $\lambda=(\lambda_1\geq\lambda_2\geq\cdots\geq \lambda_j>0)$ be a partition of $n$ with $j\geq q+1$. If
 there exists $\lambda_a\geq3$ for some $a\geq q+1$ or
 we have $\lambda_a=\lambda_{a+1}=2$ for some $a\geq q+1$,
then we have $(s \otimes \zeta_{\lambda},1)_{(Z_{\lambda})_s}=0$ for every $s\in S$.}

{\noindent\bf Proof.}
Let $$s=\prod^{d}_{i=1}(k_i,l_i)\in S$$ for $1\leq k_i<k_{i+1}\leq n-q$, $n-q+1\leq l_i <l_{i+1}\leq n$, $d\leq q$. Then we have $$s\circ\mathfrak{Q}=\{s(n-q+1),s(n-q+2),\cdots,s(n)\}.$$ We denote by $$P_{\lambda_i}=\{\lambda_1+\lambda_2+\cdots+\lambda_{i-1}+1,\lambda_1+\lambda_2+\cdots+\lambda_{i-1}+2,\cdots,\lambda_1+\lambda_2+\cdots+\lambda_{i}\}$$  the $i$-th index set of the partition $\lambda$.
\begin{enumerate}
\item Assume that there exists $\lambda_a\geq3$ for some $a\geq q+1$. Since $|s\circ\mathfrak{Q}|=q$, $a \geq q+1$, there exists some  $b\leq a$ with $\lambda_{b}\geq3$, such that $P_{\lambda_b}\cap (s\circ\mathfrak{Q}) =\emptyset$. Hence one has $sc_{\lambda_b}s^{-1}\in \mathfrak{G}$. From Lemma $2.1$, we have $(s \otimes \zeta_{\lambda},1)_{(Z_{\lambda})_s}=0$.
\item Assume that there exists $\lambda_a=\lambda_{a+1}=2$ for some $a\geq q+1$. Since $|s\circ\mathfrak{Q}|=q$, $a \geq q+1$, there exists some $b_1\leq b_2\leq a+1$ with $\lambda_{b_1}\geq\lambda_{b_2}\geq2$, such that $P_{\lambda_{b_i}}\cap (s\circ\mathfrak{Q}) =\emptyset$. If $\lambda_{b_1}=\lambda_{b_2}=2$, then
    \begin{align*}
    \nu=&
    (\lambda_{1}+\lambda_{2}+\cdots+\lambda_{b_1-1}+1,\lambda_{1}+\lambda_{2}+\cdots+\lambda_{b_2-1}+1)\\
    &(\lambda_{1}+\lambda_{2}+\cdots+\lambda_{b_1},\lambda_{1}+\lambda_{2}+\cdots+\lambda_{b_2})
    \end{align*}
    is an element of $Z_{\lambda}$ and $s\nu s^{-1}\in \mathfrak{G}$.
    Otherwise, if $\lambda_{b_1}\geq 3$, then $sc_{\lambda_{b_1}}s^{-1}\in \mathfrak{G}$. From Lemma $2.1$, we have $(s \otimes \zeta_{\lambda},1)_{(Z_{\lambda})_s}=0$.
\end{enumerate}
The lemma follows. \qb

 Let $(\mathbb{Z}/2)^n$ be the set of functions $\delta:[n]\rightarrow\mathbb{Z}/2$. Then we write $$\delta(1,2,\cdots,n)=(\delta(1),\delta(2),\cdots,\delta(n)).$$ The right $\mathfrak{S}_{n}$ action on $(\mathbb{Z}/2)^n$ is defined by
$$\delta\circ\sigma(1,2,\cdots,n)=(\delta(\sigma(1)),\delta(\sigma(2)),\cdots,\delta(\sigma(n))),$$
where $\sigma\in\mathfrak{S}_{n}$.
Then the $c_{\lambda_i}$ acts on $\delta$  by $$(\delta\circ c_{\lambda_i})(a)=\left\{\begin{array}{rl}
\delta(a+1),& \text{if }\lambda_1+\cdots+\lambda_{i-1}+1\leq a \leq \lambda_1+\cdots+\lambda_i-1,\\
\delta(a-\lambda_i+1),& \text{if }a=\lambda_1+\cdots+\lambda_i,\\
\delta(a),& \text{otherwise}
\end{array}\right.$$
for any partition $\lambda=(\lambda_1\geq\lambda_2\geq\cdots\geq \lambda_j>0)$ of $n$.
If $\lambda_i=\lambda_{i+1}$ for some $i$,  the $\nu_i$ acts on $\delta$ by $$(\delta\circ\nu_i)(a)=\left\{\begin{array}{rl}
\delta(a+\lambda_i),& \text{if }\lambda_1+\cdots+\lambda_{i-1}+1\leq a \leq \lambda_1+\cdots+\lambda_i,\\
\delta(a-\lambda_{i+1}),& \text{if }\lambda_1+\cdots+\lambda_{i}+1\leq a \leq \lambda_1+\cdots+\lambda_{i+1},\\
\delta(a),& \text{otherwise.}
\end{array}\right.$$
Now we only consider the $Z_{\lambda}$ action on $(\mathbb{Z}/2)^n$.

{\noindent\bf Definition 2.4.}
{\it We define a map  $$\{\lambda_{1}+\cdots+\lambda_{i-1}+1,\cdots,\lambda_{1}+\cdots+\lambda_{i}\}\rightarrow\mathbb{Z}/2,$$  denoted by $\delta_{\lambda_i}^{k_i}$, by $\delta_{\lambda_i}^k(a)=(\delta\circ c_{\lambda_i}^{k-1})(a)$ and call it the $i$-th cycle of $\delta$. We write $$\delta_{\lambda_1}^{k_1}\delta_{\lambda_2}^{k_2}\cdots\delta_{\lambda_j}^{k_j}=\delta\circ(c_{\lambda_1}^{k_1-1}c_{\lambda_2}^{k_2-1}\cdots c_{\lambda_j}^{k_j-1}).$$
Hence the $Z_{\lambda}$-orbit of $\delta$ is $$\{\delta_{\lambda_1}^{k_1}\delta_{\lambda_2}^{k_2}\cdots\delta_{\lambda_j}^{k_j}\circ\nu|\nu\in N_{\lambda},1\leq k_i \leq \lambda_i\}.$$
 If $\lambda_i=\lambda_{i+1}$ for some $i$, then $\nu_i\in N_{\lambda}$ acts on  $\delta_{\lambda_1}^{k_1}\delta_{\lambda_2}^{k_2}\cdots\delta_{\lambda_j}^{k_j}$ by $$\delta_{\lambda_1}^{k_1}\cdots\delta_{\lambda_i}^{k_i}\delta_{\lambda_{i+1}}^{k_{i+1}}\cdots\delta_{\lambda_j}^{k_j}\circ\nu_i
=\delta_{\lambda_1}^{k_1}\cdots\delta_{\lambda_{i}}^{k_{i+1}}\delta_{\lambda_{i+1}}^{k_{i}}\cdots\delta_{\lambda_j}^{k_j}.$$
Let
\begin{align*}B_{i}(\delta)=&\{\lambda_1+\cdots+\lambda_{i-1}+1\leq i_j\leq\lambda_{1}+\cdots+\lambda_{i}|\delta(i_j)=1\} \\
\subset &\{\lambda_{1}+\cdots+\lambda_{i-1}+1,\cdots,\lambda_{1}+\cdots+\lambda_{i}\}.
\end{align*}
We  write  $B_{i}(\delta)=\{b_1<\cdots<b_{d_i}\}$ with $d_i=|B_{i}(\delta)|$.
 We define an element in $(\mathbb{Z}_{\geq0})^{d_i}/\langle c_{d_i}\rangle$, denoted by $\chi_i(\delta)$, by
\begin{align*}
\chi_i(\delta)&=(b_2-b_1-1,\cdots, b_{d_i}-b_{d_i-1}-1, \lambda_i-b_{d_i}+b_{1})\\
&=(b_3-b_2-1,\cdots, \lambda_i-b_{d_i}+b_{1},b_2-b_1-1)\\
&=\cdots ,\end{align*} and call it the $i$-th invariant cycle of $\delta$, while $\chi_i(\delta)=(\lambda_i-1)$ if $|B_{i}(\delta)|=1$ and $\chi_i(\delta)=\emptyset$ if $B_{i}(\delta)=\emptyset$. We define a subset of $$\{\emptyset\}\sqcup(\coprod_{d}(\mathbb{Z}_{\geq0})^{d}/\langle c_{d}\rangle),$$ denoted by $\chi(\delta)$, by
$$\chi(\delta)=\{\chi_1(\delta),\chi_2(\delta),\cdots,\chi_{j}(\delta)\}$$ and call it the full invariant cycle set of $\delta$.
We define $I_{\lambda}(q)$ by
$$I_{\lambda}(q)=\{\chi(\delta)| |B_1(\delta)|+\cdots+|B_j(\delta)|=q\}.$$ }

{\bf Example 1.} Let $\lambda=(6,6,3,2,1,1)$, \begin{align*}&\delta(1,2,\cdots,19)\\
=&(1,0,1,0,0,1,0,1,0,0,1,1,0,1,0,1,1,1,0).\end{align*} Then $$\chi(\delta)=\{(1,2,0),(2,0,1),(2),(0,0),(0)\}\in I(10),$$ $\chi_1(\delta)=(1,2,0)=(0,1,2)=(2,0,1)=\chi_2(\delta)$,
 $\chi_6(\delta)=\emptyset$, $d_1=d_2=3$, $d_3=d_5=1$, $d_4=2$, $d_6=0$.\qb

 In fact, there is an injection
$$\begin{array}{cccc}i:&\mathfrak{G}\backslash \mathfrak{S}_{n}&\hookrightarrow&(\mathbb{Z}/2)^{n}\\
&s &\mapsto& \delta_s
\end{array}$$
 defined by
$$\delta_s(i)=\left\{\begin{array}{rl}0 ,& \text{if }1\leq s(i)\leq n-q,\\
1 , & \text{if }n-q+1\leq s(i) \leq n.\end{array}\right.$$
The right $Z_{\lambda}$ action on $\mathfrak{G}\backslash \mathfrak{S}_{n}$ induces an equivalence relation $"\sim"$ on $S$. In detail,  $s \sim s' \in S$ if there exists a $g \in \mathfrak{G}$ and a $z \in Z_{\lambda}$ such that $s=gs'z$. On the other hand,  the left $\mathfrak{G}$ quotient of $\mathfrak{S}_{n}\slash Z_{\lambda}$ induces a homomorphism $$\widetilde{\delta}:\mathfrak{G}\backslash \mathfrak{S}_{n}\slash Z_{\lambda}\rightarrow(\mathbb{Z}/2)^{n}\slash Z_{\lambda}$$ defined by
 \begin{align*}
 &\widetilde{\delta}((i_1,\cdots,i_{\lambda_1})(i_{\lambda_1+1},\cdots,i_{\lambda_1+\lambda_2})\cdots(i_{\lambda_1+\lambda_2+\cdots+\lambda_{j-1}+1},\cdots,i_n))\\
 =&((\delta(i_1), \cdots,\delta(i_{\lambda_1})(\delta(i_{\lambda_1+1}),\cdots,\delta(i_{\lambda_1+\lambda_2}))\cdots
 (\delta(i_{\lambda_1+\lambda_2+\cdots+\lambda_{j-1}+1}),\cdots,\delta(i_n)),
 \end{align*}
 where $$(i_1,\cdots,i_{\lambda_1})(i_{\lambda_1+1},\cdots,i_{\lambda_1+\lambda_2})\cdots(i_{\lambda_1+\lambda_2+\cdots+\lambda_{j-1}+1},\cdots,i_n)$$ is a representative of $\mathfrak{S}_{n}\slash Z_{\lambda}$ and
 $$\delta(i)=\left\{\begin{array}{rl}0, & \text{if }1\leq i\leq n-q,\\
 1, &\text{if } n-q+1\leq i \leq n. \end{array}\right.$$
 Therefore, $s \sim s' \in S$ if
 \begin{align*}
 \widetilde{\delta}(s)=\delta_s=&((\delta(s(1)), \cdots, \delta(s(\lambda_1))(\delta(s(\lambda_1+1)),\cdots,\delta(s(\lambda_1+\lambda_2))\\
 &\cdots(\delta(s(\lambda_1+\lambda_2+\cdots+\lambda_{j-1}+1)),\cdots,\delta(s(n)))\\
 =&((\delta(s'(1), \cdots,\delta(s'(\lambda_1))^{k_1}
 (\delta(s'(\lambda_1+1)),\cdots,\delta(s'(\lambda_1+\lambda_2)))^{k_2}\\
 &\cdots(\delta(s'(\lambda_1+\lambda_2+\cdots+\lambda_{j-1}+1)),\cdots,\delta(s'(n)))^{k_j})\circ\nu\\
 =&\delta_{s'}\circ(c_{\lambda_1}^{k_1-1}c_{\lambda_2}^{k_2-1}\cdots c_{\lambda_j}^{k_j-1}\nu)=\widetilde{\delta}(s').
 \end{align*}
 Thus $\chi$ is a complete invariant of $\mathfrak{G}\backslash \mathfrak{S}_{n}\slash Z_{\lambda}$, and $\widetilde{\delta}$ induces a bijection between $\mathfrak{G}\backslash \mathfrak{S}_{n}\slash Z_{\lambda}$ and $I_{\lambda}(q)$.

Consider the full invariant cycle set of an element $s\in\mathfrak{G}\backslash \mathfrak{S}_{n}\slash Z_{\lambda}$
$$\chi(\delta_{s})=\{\chi_1(\delta_{s}),\cdots,\chi_j(\delta_{s})\}.$$
 For each $\chi_{i}(\delta_{s})$, we write $\chi_i(\delta_{s})=(a_{1},a_{2},\cdots, a_{d_i})$ respectively.
Since $\chi(\delta_{s})$ is an unordered set and
\begin{align*}
\chi_{i}(\delta_{s})=&(a_{1},a_{2},\cdots,a_{d_i-1}. a_{d_i})\\
=&(a_2,a_3,\cdots,a_{d_i},a_1)\\
=&(a_3,a_4,\cdots,a_1.a_2)\\
=&\cdots,
\end{align*}
  we define a lexicographic order in $\chi(\delta_s)$ by
\begin{enumerate}
\item $\chi_i(\delta_s)=\text{min}\{(a_{1},a_{2},\cdots, a_{d_i}),(a_2,a_3,\cdots,a_{d_i},a_1), \cdots\}$, where
$$\text{min}\{(a_{1},a_{2},\cdots, a_{d_i}),(a_2,a_3,\cdots,a_{d_i},a_1), \cdots\}$$
is given by the lexicographic order. We write  $\text{min}\{\chi_i(\delta_{s}),\emptyset\}=\emptyset$.

\item  If $\lambda_i=\lambda_{i+1}$, then $|B_i(\delta_{s})|\geq |B_{i+1}(\delta_{s})|$.

\item If $\lambda_i=\lambda_{i+1}$ and $|B_i(\delta_{s})|= |B_{i+1}(\delta_{s})|$, then  $\chi_i(\delta_{s})\leq\chi_{i+1}(\delta_{s})$, where the order is the lexicographic order in $(1)$.
\end{enumerate}

{\noindent\bf Theorem 2.5.}
{\it $(s \otimes \zeta_{\lambda},1)_{(Z_{\lambda})_s}=1$ if and only if $\chi(\delta_s)$ satisfies
\begin{enumerate}
\item $1\leq |B_i(\delta_s)|\leq \lambda_i-1$ for every $\lambda_i\geq 3$.
\item If $\lambda_i\geq 3$ with $\lambda_i\equiv0,1,3$ mod $4$, then the number of minimum elements in
      $$\{(a_{1},a_{2},\cdots, a_{d_i}),(a_2,a_3,\cdots,a_{d_i},a_1),\cdots\}$$ equals $1$.
\item If $\lambda_i\geq 2$ with $\lambda_i\equiv2$ mod $4$, then the number of minimum elements in
 $$\{(a_{1},a_{2},\cdots, a_{d_i}),(a_2,a_3,\cdots,a_{d_i},a_1),\cdots\}$$ is at most $2$.
\item If there exist some $\lambda_i=\lambda_{i+1}=2k$,  then $|B_i(\delta_{s})|>|B_{i+1}(\delta_{s})|$ or $|B_i(\delta_{s})|= |B_{i+1}(\delta_{s})|$ with $\chi_i(\delta_{s})<\chi_{i+1}(\delta_{s})$.
\end{enumerate}}

{\noindent\bf Proof.}  Note that if there exist some
$$a_1+a_2+\cdots+a_{l}=k$$ such that $pk+pl=\lambda_i$ and
\begin{align*}
\chi_i(\delta_s)=&(a_1,a_2,\cdots,a_{d_i})\\
=&(a_1,a_2,\cdots,a_l,a_1,a_2,\cdots,a_l,\cdots,a_1,a_2,\cdots,a_l),
\end{align*}
then we have the number of minimum elements in
      $$\{(a_{1},a_{2},\cdots, a_{d_i}),(a_2,a_3,\cdots,a_{d_i},a_1),\cdots\}$$ is $p$ and $sc_{\lambda_i}^{\lambda_i/p}s^{-1}\in\mathfrak{G}$. By
      the definition of $\phi_{\lambda_i}$, we have $\phi_{\lambda_i}(c_{\lambda_i}^{\lambda_i/p})=1$ if and only if $p=1$ or
      $p=2$ with $\lambda_i\equiv2$ mod $4$. If there exist some $\lambda_i=\lambda_{i+1}$ such that $|B_i(\delta_{s})|= |B_{i+1}(\delta_{s})|$ with $\chi_i(\delta_{s})=\chi_{i+1}(\delta_{s})$, then we have $s\nu_i s^{-1}\in\mathfrak{G}$. By the definition of $\alpha$, we have $\alpha(\nu_i)=1$ if and only if
      $\lambda_i=\lambda_{i+1}=2k+1$.

Assume that an element $s\in \mathfrak{G}\backslash\mathfrak{S}_{n}\slash Z_{\lambda}$ satisfies $(1)$, $(2)$, $(3)$ and $(4)$. From Lemma $2.1$, we have $(s \otimes \zeta_{\lambda},1)_{(Z_{\lambda})_s}=1$.

Inversely, assume that $(s \otimes \zeta_{\lambda},1)_{(Z_{\lambda})_s}=1$.
 \begin{enumerate}
\item Assume that $\lambda_i\geq 3$. If $|B_i(\delta_s)|=0$ or $|B_i(\delta_s)|=\lambda_i$, then there exists a $c_{\lambda_i}\in C_\lambda$ such that $sc_{\lambda_i}s^{-1}\in\mathfrak{G}$. From Lemma $2.1$, we have  $(s \otimes \zeta_{\lambda},1)_{(Z_{\lambda})_s}=0$. Hence we have $1\leq |B_i(\delta_s)|\leq \lambda_i-1$ for every $\lambda_i\geq 3$.

\item  Assume that $\lambda_i\geq 3$ with $\lambda_i\equiv0,1,3$ mod $4$. If the number of minimum elements in
      $$\{(a_{1},a_{2},\cdots, a_{d_i}),(a_2,a_3,\cdots,a_{d_i},a_1),\cdots\}$$ is more than $1$, then there exists a $c_{\lambda_i}\in C_\lambda$ and a $1\leq p\leq  \lambda_i-1$ such that $sc_{\lambda_i}^{\lambda_i/p}s^{-1}\in\mathfrak{G}$.
From Lemma $2.1$, we have $(s \otimes \zeta_{\lambda},1)_{(Z_{\lambda})_s}=0$ when $\lambda_i\equiv0,1,3$ mod $4$. Hence we have the number of minimum elements in
      $$\{(a_{1},a_{2},\cdots, a_{d_i}),(a_2,a_3,\cdots,a_{d_i},a_1),\cdots\}$$ equals $1$.

\item Assume that  $\lambda_i\geq 2$ with $\lambda_i\equiv2$ mod $4$. If the number of minimum elements in
 $$\{(a_{1},a_{2},\cdots, a_{d_i}),(a_2,a_3,\cdots,a_{d_i},a_1),\cdots\}$$ is more than $2$, then similar with the case $(2)$, from Lemma $2.1$, we have $(s \otimes \zeta_{\lambda},1)_{(Z_{\lambda})_s}=0$. We remark that if the number of minimum elements in
 $$\{(a_{1},a_{2},\cdots, a_{d_i}),(a_2,a_3,\cdots,a_{d_i},a_1),\cdots\}$$ equals  $2$,
  then  there exists a $c_{\lambda_i}\in C_\lambda$  such that $sc_{\lambda_i}^{\lambda_i/2}s^{-1}\in\mathfrak{G}$. Since
    $$\phi_{\lambda_i}(c_{\lambda_i}^{\lambda_i/2})=-e^{\pi i}=1,$$ we can get $(s \otimes \zeta_{\lambda},1)_{(Z_{\lambda})_s}=1$ for some $s\in\mathfrak{G}\backslash \mathfrak{S}_{p+q}\slash Z_{\lambda}$. Hence the number of minimum elements in
 $$\{(a_{1},a_{2},\cdots, a_{d_i}),(a_2,a_3,\cdots,a_{d_i},a_1),\cdots\}$$ is at most $2$.

\item Assume that there exist some $\lambda_i=\lambda_{i+1}=2k$ such that $|B_i(\delta_{s})|= |B_{i+1}(\delta_{s})|$ with $\chi_i(s)=\chi_{i+1}(s)$, then there exists a $\nu_i\in N_\lambda$ such that
$s\nu_i s^{-1}\in\mathfrak{G}$. From Lemma $2.1$, we have $(s \otimes \zeta_{\lambda},1)_{(Z_{\lambda})_s}=0$.
Hence we have $|B_i(\delta_{s})|>|B_{i+1}(\delta_{s})|$ or $|B_i(\delta_{s})|= |B_{i+1}(\delta_{s})|$ with $\chi_i(s)<\chi_{i+1}(s)$.
\end{enumerate}

The theorem follows.
 \qb

Let $$H^{\lambda_1+\lambda_2+\lambda_a-a}(P_n)^\mathfrak{G}|_{(\lambda_1,\lambda_2,\cdots,\lambda_a)}\cong\bigoplus_{(s \otimes \zeta_{\lambda},1)_{(Z_{\lambda})_s}=1 \atop s\in \mathfrak{G}\backslash\mathfrak{S}_n / Z_{\lambda}}(s \otimes \zeta_{\lambda})$$ be the subgroup of $$H^{\lambda_1+\lambda_2+\lambda_a-a}(P_n)^\mathfrak{G}\cong\bigoplus_{\lambda}\bigoplus_{(s \otimes \zeta_{\lambda},1)_{(Z_{\lambda})_s}=1 \atop s\in \mathfrak{G}\backslash\mathfrak{S}_n / Z_{\lambda}}(s \otimes \zeta_{\lambda})$$ restricted by the partition $\lambda=(\lambda_1,\lambda_2,\cdots,\lambda_a,1,\cdots)$ for some $a\geq1$, where $\lambda_1 \geq \lambda_2\geq \cdots\geq \lambda_a\geq 2$.
It is easy to see that if $\lambda=(\lambda_1,1,1,\cdots)$, then $$H^{\lambda_1-1}(P_n)^\mathfrak{G}|_{(\lambda_1)}\cong\bigoplus_{(s \otimes \phi_{\lambda},1)_{(Z_{\lambda})_s}=1\atop s\in \mathfrak{G}\backslash\mathfrak{S}_n / Z_{\lambda}}(s \otimes \phi_{\lambda_1}).$$
We write the partitions $\widetilde{\lambda}_i=(\lambda_i,1,1,\cdots)$ and
$$\widetilde{N}_\lambda=\langle\nu_i\in N_{\lambda}|\lambda_i=\lambda_{i+1}\geq2\rangle$$ the subgroup of $N_\lambda$ .  The $\widetilde{N}_\lambda$ action on $I_{\widetilde{\lambda}_1}(q)\times I_{\widetilde{\lambda}_2}(q)\cdots \times I_{\widetilde{\lambda}_i}(q)$ is given by
$$
\nu_i\circ(\cdots,\chi(\delta_{s_i}),\chi(\delta_{s_{i+1}}),\cdots)
=(\cdots,\chi(\delta_{s_{i+1}}),\chi(\delta_{s_{i}}),\cdots)
$$
for some $\lambda_i=\lambda_{i+1}\geq2$, $\nu_i\in N_{\lambda}$,
while it induces the $\widetilde{N}_\lambda$ action on $$\bigotimes_{i=1}^{a}H^{\lambda_i-1}(P_n)^\mathfrak{G}|_{(\lambda_i)}$$ by
$$
\nu_i\circ(\cdots\otimes(s_i\otimes\phi_{\lambda_i})\otimes(s_{i+1}\otimes\phi_{\lambda_{i+1}})\otimes\cdots)
=\cdots\otimes(s_{i+1}\otimes\phi_{\lambda_{i}})\otimes(s_{i}\otimes\phi_{\lambda_{i+1}})\otimes\cdots.
$$

{\noindent\bf Theorem 2.6.}
{\it Let $\lambda=(\lambda_1,\lambda_2,\cdots,\lambda_a,1,\cdots)$. There is a surjection of $\mathbb{Q}$-modules
\[  \varphi: \bigotimes_{i=1}^{a}H^{\lambda_i-1}(P_n)^\mathfrak{G}|_{(\lambda_i)}
\rightarrow H^{\lambda_1+\lambda_2+\cdots+\lambda_a-a}(P_n)^\mathfrak{G}|_{(\lambda_1,\lambda_2,\cdots,\lambda_a)},\]
defined by \begin{align*}
&\varphi((s_1\otimes\phi_{\lambda_1})\otimes\cdots\otimes(s_a\otimes\phi_{\lambda_a}))\\
=&\left\{\begin{array}{cl}
s_0\otimes(\alpha_{\lambda}\phi_{\lambda_1}\cdots\phi_{\lambda_a}), &\text{ if } \sum_{i=1}^{a}(\lambda_i)-n+q\leq\sum_{i=1}^a|B_1(\delta_{s_i})|\leq q,\\
0, & \text{ otherwise, }\end{array}\right.
\end{align*}
where $s_0$ is induced by $\chi_i(\delta_{s_0})=\chi_1(\delta_{s_i})$.
 The kernel of $\varphi$ is generated by the union of
\begin{align*}
A_1=&\{\bigotimes_{i=1}^a(s_i\otimes\phi_{\lambda_i})-\nu\circ\bigotimes_{i=1}^a(s_i\otimes\phi_{\lambda_i})|\nu\in\widetilde{N}_\lambda \},\\
A_2=&\{\bigotimes_{i=1}^a(s_i\otimes\phi_{\lambda_i})| \chi_1(\delta_{s_{i}})=\chi_1(\delta_{s_{j}}),\text{ if }\lambda_i=\lambda_j=2k\},\\
A_3=&\{\bigotimes_{i=1}^a(s_i\otimes\phi_{\lambda_i})|\sum_{i=1}^a|B_1(\delta_{s_i})|<\sum_{i=1}^{a}(\lambda_i)-n+q\text{ or } \sum_{i=1}^a|B_1(\delta_{s_i})|>q \}.
\end{align*}
The image of a set of representatives of $\widetilde{N}_{\lambda}$ action on
$$A_4=\left\{\bigotimes_{i=1}^{a}(s_{i}\otimes\phi_{\lambda_i})\left|\begin{array}{c}\chi_1(\delta_{s_{i}})\neq\chi_1(\delta_{s_{j}}),\text{ if }\lambda_i=\lambda_j=2k, \\ \sum_{i=1}^{a}(\lambda_i)-n+q\leq\sum_{i=1}^a|B_1(\delta_{s_i})|\leq
q\end{array}\right.\right\}$$
is a basis of
$H^{\lambda_1+\lambda_2+\cdots+\lambda_a-a}(P_n)^\mathfrak{G}|_{(\lambda_1,\lambda_2,\cdots,\lambda_a)}$.
}

{\noindent\bf Proof.} Note that $\bigotimes_{i=1}^{a}H^{\lambda_i-1}(P_n)^\mathfrak{G}|_{(\lambda_i)}$ and $H^{\lambda_1+\lambda_2+\cdots+\lambda_a-a}(P_n)^\mathfrak{G}|_{(\lambda_1,\lambda_2,\cdots,\lambda_a)}$ are finitely generated as $\mathbb{Q}$-vector spaces.
$$
\{\bigotimes_{i=1}^{a}(s_{i}\otimes\phi_{\lambda_i})|(s_{i}\otimes\phi_{\lambda_i},1)_{(Z_{\lambda_i})_{s_i}}=1\}
$$ is  a set of  generators of $\bigotimes_{i=1}^{a}H^{\lambda_i-1}(P_n)^\mathfrak{G}|_{(\lambda_i)}$ and
$$\{s_{0}\otimes(\alpha_{\lambda}\phi_{\lambda_1}\cdots\phi_{\lambda_a})|(s_{0}\otimes(\alpha_{\lambda}\phi_{\lambda_1}\cdots\phi_{\lambda_a}),1)_{(Z_{\lambda})_{s_0}}=1\}$$
is a set of  generators of $H^{\lambda_1+\lambda_2+\cdots+\lambda_a-a}(P_n)^\mathfrak{G}|_{(\lambda_1,\lambda_2,\cdots,\lambda_a)}$. For every generator $s_0\otimes(\alpha_{\lambda}\phi_{\lambda_1}\cdots\phi_{\lambda_a})$, consider the full invariant cycle set of $s_0$.
Note that for a  $$\chi(\delta_{s_0})=\{\chi_1(\delta_{s_0}),\cdots,\chi_a(\delta_{s_0}),\cdots\}\in I_{\lambda}(q),$$ there exists a unique
 $$(\chi(\delta_{s_1}),\chi(\delta_{s_2}),\cdots,\chi(\delta_{s_a}))\in(I_{\widetilde{\lambda}_1}(q)\times I_{\widetilde{\lambda}_2}(q)\cdots \times I_{\widetilde{\lambda}_i}(q))$$ such that
$$\chi(\delta_{s_i})=\{\chi_1(\delta_{s_i})=\chi_i(\delta_{s_0}),(0),\cdots,(0)\},$$
 where the number of $(0)$ is $q-|B_1(s_i)|$.
Hence the map
$$
\widetilde{\varphi}: H^{\lambda_1+\lambda_2+\cdots+\lambda_a-a}(P_n)^\mathfrak{G}|_{(\lambda_1,\lambda_2,\cdots,\lambda_a)}\rightarrow \bigotimes_{i=1}^{a}H^{\lambda_i-1}(P_n)^\mathfrak{G}|_{(\lambda_i)}
 $$
defined by  $$\widetilde{\varphi}((s_{0}\otimes(\alpha_{\lambda}\phi_{\lambda_1}\cdots\phi_{\lambda_a})))=\bigotimes_{i=1}^{a}(s_{i}\otimes\phi_{\lambda_i}),$$
 where $\chi_1(\delta_{s_{i}})=\chi_i(\delta_{s_{0}})$, $1\leq i \leq a$, is well-defined.
   By the definition of $\widetilde{\varphi}$,
 we have $\varphi\widetilde{\varphi}=id$. Hence we have $\varphi$ is surjective and $\widetilde{\varphi}$ is injective.

Now we consider the set of generators of $H^{\lambda_1+\lambda_2+\cdots+\lambda_a-a}(P_n)^\mathfrak{G}|_{(\lambda_1,\lambda_2,\cdots,\lambda_a)}$.
 Note that for every
$$\widetilde{\varphi}(s_{0}\otimes(\alpha_{\lambda}\phi_{\lambda_1}\cdots\phi_{\lambda_a}))=\bigotimes_{i=1}^{a}(s_{i}\otimes\phi_{\lambda_i}),$$
we have $(s_{i}\otimes\phi_{\lambda_i},1)_{(Z_{\lambda_i})_{s_i}}=1$ and $$\sum_{i=1}^{a}(\lambda_i)-n+q\leq\sum_{i=1}^a|B_1(\delta_{s_i})|\leq
q.$$
If $\lambda_i=\lambda_j=2k$ for some $i$ and $j$, we also have $\chi_1(\delta_{s_{i}})\neq\chi_1(\delta_{s_{j}})$.
Hence  the image of
$$A_4=\left\{\bigotimes_{i=1}^{a}(s_{i}\otimes\phi_{\lambda_i})\left|\begin{array}{c}\chi_1(\delta_{s_{i}})\neq\chi_1(\delta_{s_{j}}),\text{ if }\lambda_i=\lambda_j=2k, \\ \sum_{i=1}^{a}(\lambda_i)-n+q\leq\sum_{i=1}^a|B_1(\delta_{s_i})|\leq
q\end{array}\right.\right\}$$
 generates  $H^{\lambda_1+\lambda_2+\cdots+\lambda_a-a}(P_n)^\mathfrak{G}|_{(\lambda_1,\lambda_2,\cdots,\lambda_a)}$.
 We write $\widetilde{A}_4$ a set of representatives of $\widetilde{N}_\lambda$ action on $A_4$.
Since $\chi(\delta_{s_0})=\chi(\delta_{s_0'})$ if there exists a $\nu\in\widetilde{N}_\lambda$ such that $s_0=\nu\circ s_0'$, we have $\widetilde{A}_4$ is a basis of $im\widetilde{\varphi}$ and the image of $\widetilde{A}_4$ is a basis of $H^{\lambda_1+\lambda_2+\cdots+\lambda_a-a}(P_n)^\mathfrak{G}|_{(\lambda_1,\lambda_2,\cdots,\lambda_a)}$.

Finally we consider the kernel of $\varphi$. Since $\widetilde{\varphi}$ is injective, we have
$$\bigotimes_{i=1}^{a}H^{\lambda_i-1}(P_n)^\mathfrak{G}|_{(\lambda_i)}=ker\varphi \oplus im\widetilde{\varphi}.$$
Note that the intersection of $\mathbb{Q}A_1$ and $\mathbb{Q}\widetilde{A}_4$ is $0$, $\mathbb{Q}A_1\oplus\mathbb{Q}\widetilde{A}_4=\mathbb{Q}A_4$, $A_i\cap A_4=\emptyset$, $2\leq i \leq3$, and
$$ A_2\cup A_3\cup A_4=\{\bigotimes_{i=1}^{a}(s_{i}\otimes\phi_{\lambda_i})|(s_{i}\otimes\phi_{\lambda_i},1)_{(Z_{\lambda_i})_{s_i}}=1\}.$$
Thus we have $A_1 \cup A_2\cup A_3$ is a set of generators of $ker\varphi$.

 The theorem follows.\qb

{\bf Remark:}
Suppose that $\lambda_1+\cdots+\lambda_a\leq n-q$ and $|B_1(\widetilde{\delta}(s_1))|+\cdots+ |B_a(\widetilde{\delta}(s_a))|\leq q$. Since we have $$\chi_1(\widetilde{\delta}(\prod^{d_j}_{i=1}(k_{ji},l_{ji})))=\chi_i (\widetilde{\delta}(\prod^{d_j}_{i=1}(\sum_{b=1}^{j-1}\lambda_j+k_{bi},n-q+\sum_{b=1}^{j-1}d_b+i)))$$ for every $d_j$, the surjection $\varphi$ can be determined by \begin{align*}
&\varphi((\prod^{d_1}_{i=1}(k_{1i},l_{1i})\otimes\phi_{\lambda_1})\otimes\cdots\otimes(\prod^{d_a}_{i=1}(k_{ai},l_{ai})\otimes\phi_{\lambda_a}))\\
=&(\prod^{d_1}_{i=1}(k_{1i},n-q+i)\cdots\prod^{d_a}_{i=1}(\sum_{j=1}^{a-1}\lambda_j+k_{ai},n-q+\sum_{j=1}^{a-1}d_j+i))\otimes \zeta_{\lambda},
\end{align*}
where $\zeta_{\lambda}=\alpha_{\lambda}\phi_{\lambda_1}\phi_{\lambda_2}\cdots\phi_{\lambda_a}$, $1\leq k_{ji}\leq \lambda_j$, $n-q+1\leq l_{ji}\leq p+q$.

{\noindent\bf Theorem 2.7.}
{\it Let $\mathfrak{G}_1=\mathfrak{S}_{n-q}\times\mathfrak{S}_{q}$ and $\mathfrak{G}_2=\mathfrak{S}_{n-q-1}\times\mathfrak{S}_{q+1}$, $n-q-1\geq q+1$.
In the range $*\leq n-q-2$, there exists an injection
\[  \psi: H^{*}(P_n)^{\mathfrak{G}_{1}}
\rightarrow H^{*}(P_n)^{\mathfrak{G}_{2}},\]
where $\psi( s \otimes \zeta_{\lambda})=s' \otimes \zeta_{\lambda}$ is given by $\chi(\delta_{s'})=\{\chi(\delta_{s}),(0)\}$.
Furthermore, $\psi$ is a bijection in the range $*\leq q-1$.}

{\noindent\bf Proof.} Since $H^{*}(P_n)^{\mathfrak{G}_1}$ is finitely generated as a $\mathbb{Q}$-vector space,
$$\bigcup_{\lambda}\{s\otimes \zeta_{\lambda}|(s \otimes \zeta_{\lambda},1)_{(Z_{\lambda})_s}=1\}$$
is a set of generators of $H^{*}(P_n)^{\mathfrak{G}_1}$.

  In order to define the map $\psi$, we will show that $\lambda_j=1$ and $\chi_j(\delta_s)=\emptyset$ for every partition  $\lambda=(\lambda_1\geq\lambda_2\geq\cdots\geq\lambda_j)$ with $(s \otimes \zeta_{\lambda},1)_{(Z_{\lambda})_s}=1$ in the range $*=n-j\leq n-q-2$.

  Consider a full invariant cycle set
$$\chi(\delta_s)=\{\chi_1(\delta_s),\chi_2(\delta_s),\cdots, \chi_j(\delta_s)\}\in I_{\lambda}(q)$$
 associated with a generator $s\otimes \zeta_{\lambda}$.
   From Lemma $2.3$, we have $\lambda_j=1$ when $j\geq q+2$. Since $j\geq q+2$, there exist at least two $\lambda_i$, $1\leq i\leq j$, such that
$\chi_i(\delta_s)=\emptyset$.  From Theorem $2.5$, when $\lambda_i\geq3$, we have $\chi_i(\delta_s)\neq\emptyset$. When $\lambda_i=2$, we have the number of $\chi_i(\delta_s)=\emptyset$  is at most one. Thus we have $\chi_j(\delta_s)=\emptyset$. Hence $\psi$ is well-defined.

Note that $\{\chi(\delta_{s}),(0)\}$ is uniquely determined by $\chi(\delta_{s})$, while $(s' \otimes \zeta_{\lambda},1)_{(Z_{\lambda})_{s'}}=(s \otimes \zeta_{\lambda},1)_{(Z_{\lambda})_{s}}$ for every $\chi(\delta_{s'})=\{\chi(\delta_{s}),(0)\}$. Hence we have $\psi$ is injective in the range $*=n-j\leq n-q-2$.

When $*\geq q$, consider the partition $\lambda=(\lambda_1,2,1,\cdots)$, one has $\lambda_1\geq q$. Note that there exists a   full invariant cycle set $$\chi(\delta_s)=\{\chi_1(\delta_s)=(0,0,\cdots,0,\lambda_1-q+1),\chi_2(\delta_s)=(0,0)\}\in I_{\lambda}(q+1)$$
  such that $(s \otimes \zeta_{\lambda},1)_{(Z_{\lambda})_{s}}=1$ and $|B_1(\delta_s)|+|B_2(\delta_s)|=q+1$. It is easy to see that $s \otimes \zeta_{\lambda}$ is not a image of $\psi$, hence $\psi$ is not surjective when $*\geq q$.

Now assume that $*\leq q-1$, from Theorem $2.5$,
we need to consider the partitions
\begin{enumerate}
\item[(i)] $\lambda=(\lambda_1,\lambda_2,\cdots, \lambda_{k},1,1,\cdots)$, $\lambda_i\geq3$, $\lambda_1+\lambda_2+\cdots+\lambda_k\leq q+k-1$.

\item[(ii)] $\lambda=(\lambda_1,\lambda_2,\cdots, \lambda_{k},2,1,1,\cdots)$, $\lambda_i\geq3$, $\lambda_1+\lambda_2+\cdots+\lambda_k\leq q+k-2$.

\item[(iii)] $\lambda=(\lambda_1,\lambda_2,\cdots, \lambda_{k},2,2,1,1,\cdots)$, $\lambda_i\geq3$, $\lambda_1+\lambda_2+\cdots+\lambda_k\leq q+k-3$.

\item[(iv)] $\lambda=(\lambda_1,\lambda_2,\cdots, \lambda_{k},2,2,2,1,1,\cdots)$, $\lambda_i\geq3$, $\lambda_1+\lambda_2+\cdots+\lambda_k\leq q+k-4$.
\end{enumerate}

\begin{enumerate}
\item For case (i), from Theorem $2.5$, one has $$\sum_{i=1}^k|B_i(\delta_s)|\leq \sum_{i=1}^{k}(\lambda_i-1)\leq q-1$$ for every $(s \otimes \zeta_{\lambda},1)_{(Z_{\lambda})_{s}}=1$.

\item For case (ii), from Theorem $2.5$ and Theorem $2.6$, one has
$$\sum_{i=1}^{k+1}|B_i(\delta_s)|\leq \sum_{i=1}^{k}(\lambda_i-1)+2\leq q$$  for every $(s \otimes \zeta_{\lambda},1)_{(Z_{\lambda})_{s}}=1$.
We remark that $\chi_{k+1}(\delta_s)$ can be $(0,0)$ such that $(s \otimes \zeta_{\lambda},1)_{(Z_{\lambda})_{s}}=1$.

\item For case (iii), from Theorem $2.5$ and Theorem $2.6$, one has
$$\sum_{i=1}^{k+2}|B_i(\delta_s)|\leq \sum_{i=1}^{k}(\lambda_i-1)+3\leq q$$  for every $(s \otimes \zeta_{\lambda},1)_{(Z_{\lambda})_{s}}=1$.
We remark that $$\{\chi_{k+1}(\delta_s),\chi_{k+2}(\delta_s)\}$$ can be $\{(0,0),(1)\}$ such that $(s \otimes \zeta_{\lambda},1)_{(Z_{\lambda})_{s}}=1$.

\item For case (iv), from Theorem $2.5$ and Theorem $2.6$, one has
$$\sum_{i=1}^{k+3}|B_i(\delta_s)|\leq \sum_{i=1}^{k}(\lambda_i-1)+3\leq q-1$$  for every $(s \otimes \zeta_{\lambda},1)_{(Z_{\lambda})_{s}}=1$.
We remark that $$\{\chi_{k+1}(\delta_s),\chi_{k+2}(\delta_s),\chi_{k+3}(\delta_s)\}$$ can be $\{(0,0),(1),\emptyset\}$ such that $(s \otimes \zeta_{\lambda},1)_{(Z_{\lambda})_{s}}=1$.
\end{enumerate}

Now we write $\lambda=(\lambda_1,\lambda_2,\cdots, \lambda_{k},1,1,\cdots)$, $\lambda_i\geq 2$, with respect to those cases, then one has
$$\sum_{i=1}^{k}|B_i(\delta_s)|\leq q$$ for every $(s \otimes \zeta_{\lambda},1)_{(Z_{\lambda})_{s}}=1$.
Since $\chi(\delta_s)\in I_{\lambda}(q+1)$, there exists at least one $\chi_{i}(\delta_s)=(0)$ for every $(s \otimes \zeta_{\lambda},1)_{(Z_{\lambda})_{s}}=1$.
Hence there exists a $\chi(\delta_{s'})\in I_{\lambda}(q)$ such that $\{\chi(\delta_{s'}),(0)\}=\chi(\delta_{s})$. Thus $\psi$ is surjective when $*\leq q-1$.

Since $q+1\leq n-q-1$, we complete our proof.\qb

\section{ the dimension of  $H^{*}(P_n)^\mathfrak{G}$}

Recall that $\lambda=(\lambda_1,\lambda_2,\cdots,\lambda_a,1,1,\cdots)$, $\lambda_i\geq 2$, $1\leq i \leq a$, is a $(n-*)$ partition of $n$. We denote by $\Lambda(n-*)$ the set of all  $(n-*)$ partitions of $n$.
For a partition $\lambda\in\Lambda(n-*)$, we consider $$d=(d_1,d_2,\cdots,d_a)\in(\mathbb{Z}_{\geq0})^a$$ and
$$(\lambda,d)=((\lambda_1,d_1),(\lambda_2,d_2),\cdots,(\lambda_a,d_a))$$ which satisfies $d_i\geq d_{i+1}$ if $\lambda_i=\lambda_{i+1}$. If $d$ satisfies
$d_1+d_2+\cdots+d_a\leq q$, then it is easy to see that there exists some full invariant cycle sets $\chi(\delta_s)$ such that $|B_i(\delta_s)|=d_i$.
In order to calculate the dimension of $H^{*}(P_n)^\mathfrak{G}$,  we define
$$\Pi(\lambda_i,d_i)=\{\chi_i(\delta_s)\in(\mathbb{Z}_{\geq0})^{d_i}/\langle c_{d_i}\rangle||B_i(\delta_s)|=d_i,(s\otimes\phi_{\lambda_i},1)_{(Z_{\lambda})_s}=1\}$$
and denote by $|\Pi(\lambda_i,d_i)|$ the cardinality of $\Pi(\lambda_i,d_i)$.
We denote
\begin{align*}
M&=\{(\lambda,d)|\lambda=(\lambda_1,\lambda_2,\cdots,\lambda_a,1,1,\cdots)\in\Lambda(n-*),\sum_{i=1}^{a}d_i\leq q\},\\
M_1(\lambda,d)&=\{1\leq i \leq a|  \lambda_i\equiv 1 \text{ mod } 2\},\\
M_2(\lambda,d)&=\{1\leq i \leq a|\lambda_i \equiv 0 \text{ mod } 2\},\\
|(\lambda_i,d_i)|&=\sharp\{1\leq j\leq a|(\lambda_j,d_j)=(\lambda_{i},d_i)\}.
\end{align*}
For $1\leq b \leq |(\lambda_i,d_i)|$, we also define simply
$$
K(b)=\{(k_1,k_2,\cdots,k_b)\in (\mathbb{Z}_{>0})^b| k_1+k_2+\cdots+k_b=|(\lambda_i,d_i)|\}
$$
for $K(\lambda_i,d_i,b)$.

{\noindent\bf Theorem 3.1.}
{\it The dimension of $H^{*}(P_n)^\mathfrak{G}$ equals
$$
\sum_{(\lambda,d)\in M}\prod_{i\in M_1(\lambda,d)}
(\sum_{b=1}^{|(\lambda_i,d_i)|}|K( b)|\left(\begin{array}{c}|\Pi(\lambda_i,d_i)|\\b\end{array}\right))
\prod_{i\in M_2(\lambda,d)}\left(\begin{array}{c}|\Pi(\lambda_i,d_i)|\\|(\lambda_i,d_i)|\end{array}\right),
$$
where
$\left(\begin{array}{c}|\Pi(\lambda_i,d_i)|\\b\end{array}\right)$, $1\leq b\leq |\Pi(\lambda_i,d_i)|$, is the binomial coefficient and
 we write
$\left(\begin{array}{c}|\Pi(\lambda_i,d_i)|\\b\end{array}\right)=0$
if $|\Pi(\lambda_i,d_i)|< b$.}

{\noindent\bf Proof.} By a direct application of Theorem $2.5$ and Theorem $2.6$, we only need to consider $(\lambda,d)\in M$. For a fixed $(\lambda,d)$, we only need to consider each of
$$(\lambda_i,d_i)\in\{(\lambda_1,d_1),(\lambda_2,d_2),\cdots,(\lambda_a,d_a)\}.$$

If $i\in M_1(\lambda,d)$, by the definition of the lexicographic order in $\chi(\delta_s)$, we assume that
\begin{align*}
\chi_{i}(\delta_s)&=\chi_{i+1}(\delta_s)=\cdots=\chi_{i+k_1-1}(\delta_s),\\
\chi_{i+k_1}(\delta_s)&=\chi_{i+k_1+1}(\delta_s)=\cdots=\chi_{i+k_1+k_2-1}(\delta_s),\\
&\cdots\\
\chi_{i+k_1+\cdots+k_{b-1}}(\delta_s)&=\chi_{i+k_1+\cdots+k_{b-1}+1}(\delta_s)=\cdots=\chi_{i+|(\lambda_i,d_i)|-1}(\delta_s),
\end{align*}
where $$\chi_{i}(\delta_s)<\chi_{i+k_1}(\delta_s)<\cdots<\chi_{i+k_1+\cdots+k_{b-1}}\in \Pi(\lambda_i,d_i).$$
Then we have $(k_1,k_2,\cdots,k_b)\in K( b)$ and the number of those elements stated above is $\left(\begin{array}{c}|\Pi(\lambda_i,d_i)|\\b\end{array}\right)$.

 If $i\in M_2(\lambda,d)$, then from Theorem $2.5$, $(4)$, we have $$\chi_i(\delta_s)<\chi_{i+1}(\delta_s)<\cdots<\chi_{i+|(\lambda_i,d_i)|-1}(\delta_s) \in\Pi(\lambda_i,d_i).$$ Hence we have the number of those elements stated above is $\left(\begin{array}{c}|\Pi(\lambda_i,d_i)|\\|(\lambda_i,d_i)|\end{array}\right)$.

Since $1\leq b\leq |\Pi(\lambda_i,d_i)|$, $(k_1,k_2,\cdots,k_b)$ runs over $K(b)$, $(\lambda,d)$ runs over $M$, we complete our proof. \qb

{\bf Remark:} For every $(\lambda,d)\in M$, we have $0\leq d_1+\cdots+d_a\leq q$. Furthermore, since $|\Pi(\lambda_i,d_i)|$ or $$\left(\begin{array}{c}|\Pi(\lambda_i,d_i)|\\b\end{array}\right)$$ can be $0$, the range of $d_1+\cdots+d_{a}$ can be smaller for some $(\lambda,d)\in M$.

 We write $\Pi(\lambda,d)$ respectively. From Theorem $2.5$, we have $\Pi(\lambda,0)=\Pi(\lambda,\lambda)=0$ for every $\lambda\geq3$ and $\Pi(2,2)=\Pi(2,1)=\Pi(2,0)=1$.

 Now we suppose $\lambda\geq3$, $1\leq d\leq \lambda-1$. For $1\leq p\leq d$, we consider
\begin{align*}
t&=(t_1,t_2,\cdots,t_p)\in(\mathbb{Z}_{>0})^p,\\
 a&=(a_{i_1}<a_{i_2}<\cdots<a_{i_p})\in(\mathbb{Z}_{\geq0})^p
\end{align*}
and the pair $T=(t,a)$.
 It is easy to see that if the pair $T$ satisfies $$t_1+t_2+\cdots+t_p=d$$ and
 $$t_1a_{i_1}+t_2a_{i_2}+\cdots+t_pa_{i_p}=\lambda-d,$$ then there exist some invariant cycles $(a_{1},\cdots,a_{d})$ such that
 $$\{a_1,a_2,\cdots,a_d\}=\{a_{i_1}<a_{i_2}<\cdots<a_{i_p}\}.$$
 We write $T(p)$ the set of all $T$ which satisfy the above conditions.

   If we calculate the number of all
 $\chi_i(\delta_s)$ which satisfy the conditions (1), (2) and (3) in Theoem $2.5$ associated with fixed  $\lambda\geq 3$, $1\leq d\leq \lambda-1$ and $T\in T(p)$,  then we obtain $|\Pi(\lambda,d)|$ by considering all $T\in T(p)$ for all $1\leq p \leq d$.

For  $T=(t,a)\in T(p)$ and a divisor $k \geq 1$ of $gcd(t_1,t_2,\cdots,t_p)$, we write $$gcd(t_1,t_2,\cdots,t_p)=kp_1^{l_1}p_2^{l_2}\cdots p_m^{l_m},$$
where $p_j$, $1\leq j\leq m$, are distinct prime numbers.
Then we denote $$P(T,k)=\{p_1,p_2,\cdots,p_m\}.$$

{\noindent\bf Definition 3.2.}
{\it For $0 \leq j \leq m$, we denote by $P(T,k)_j$ the set consisting of all subsets of $P(T,k)$ with cardinality $j$.
We agree $P(T,k)_0=\emptyset$.
For a  divisor $k\geq 1$ of $gcd(t_1,t_2,\cdots,t_p)$ and $0\leq j \leq m$, we define the number $\pi_k(T,j)$ by
\begin{align*}\pi_k(T,j)&=\sum_{\{p_{x_1},p_{x_2},\cdots, p_{x_j}\}\in P(T,k)_{j}}\left(\begin{array}{c}
\frac{d}{kp_{x_1}p_{x_2}\cdots p_{x_j}}\\
\frac{t_1}{kp_{x_1}p_{x_2}\cdots p_{x_j}},
\frac{t_2}{kp_{x_1}p_{x_2}\cdots p_{x_j}},
\cdots,
\frac{t_p}{kp_{x_1}p_{x_2}\cdots p_{x_j}}
\end{array}\right),
\end{align*}
while
\begin{align*}
\pi_k(T,0)&=\left(\begin{array}{c}
\frac{d}{k}\\
\frac{t_1}{k},
\frac{t_2}{k},
\cdots,
\frac{t_p}{k}
\end{array}\right).
\end{align*}
Here $\left(\begin{array}{c}
\frac{d}{kp_{x_1}p_{x_2}\cdots p_{x_j}}\\
\frac{t_1}{kp_{x_1}p_{x_2}\cdots p_{x_j}},
\frac{t_2}{kp_{x_1}p_{x_2}\cdots p_{x_j}},
\cdots,
\frac{t_p}{kp_{x_1}p_{x_2}\cdots p_{x_j}}
\end{array}\right)$ and $\left(\begin{array}{c}
\frac{d}{k}\\
\frac{t_1}{k},
\frac{t_2}{k},
\cdots,
\frac{t_p}{k}
\end{array}\right)$ are the multinomial coefficients.
}

{\noindent\bf  Theorem 3.3.}
{\it
\begin{enumerate}
\item[\bf 3.3.1.] If  $\lambda\equiv 0,1,3$ mod $4$, then $$|\Pi(\lambda,d)|=\sum_{p=1}^{d}\sum_{T\in T(p)}\frac{1}{d}(\sum_{j=0}^{m}(-1)^j\pi_1(T,j)).$$

\item[\bf 3.3.2.] If  $\lambda\equiv2$ mod $4$, then
$$
|\Pi(\lambda,d)|=\sum_{p=1}^{d}\sum_{T\in T(p)}\frac{1}{d}(\sum_{j=0}^{m}(-1)^j\pi_1(T,j))
+\sum_{p=1}^{d}\sum_{T\in T(p)}\frac{2}{d}(\sum_{j=0}^{m}(-1)^j\pi_2(T,j)).
$$
\end{enumerate}
}

{\noindent\bf Proof.}  In order to prove this Theorem, it suffices to show that for every
$$gcd(t_1,t_2,\cdots,t_p)=kp_1^{l_1}p_2^{l_2}\cdots p_m^{l_m}$$
with $t_1+t_2+\cdots+t_p=d$ and $t_1a_{i_1}+t_2a_{i_2}+\cdots+t_pa_{i_p}=\lambda-d$,
$$\sum_{j=0}^{m}(-1)^j\pi_k(T,j)$$
is the number of those elements  $$(a_1,a_2,\cdots,a_d)\in(\mathbb{Z}_{\geq0})^{d}$$ such that the number of minimum elements in
 $$\{(a_{1},a_{2},\cdots, a_{d}),(a_2,a_3,\cdots,a_{d},a_1),\cdots\}$$ equals $k$.

For a fixed $T=(t,a)\in T(p)$, let $(a_{1},a_{2},\cdots, a_{d})\in(\mathbb{Z}_{\geq0})^{d}$ such that $$\{a_1,a_2,\cdots,a_d\}=\{a_{i_1}<a_{i_2}<\cdots<a_{i_p}\},$$
 $$t_1+t_2+\cdots+t_p=d$$ and $$t_1a_{i_1}+t_2a_{i_2}+\cdots+t_pa_{i_p}=\lambda-d.$$ Then the number of minimum elements in
 $$\{(a_{1},a_{2},\cdots, a_{d}),(a_2,a_3,\cdots,a_{d},a_1),\cdots\}$$  is a common divisor of $gcd(t_1,t_2,\cdots,t_p)$.

 We denote by $\mathfrak{O}_0$ the set of those elements  $(a_1,a_2,\cdots,a_d)\in(\mathbb{Z}_{\geq0})^{d}$ such that the number of minimum elements in
 $$\{(a_{1},a_{2},\cdots, a_{d}),(a_2,a_3,\cdots,a_{d},a_1),\cdots\}$$ equals $k$,
 $\mathfrak{O}_i$ the set of those elements  $(a_1,a_2,\cdots,a_d)\in(\mathbb{Z}_{\geq0})^{d}$ such that the number of minimum elements in
 $$\{(a_{1},a_{2},\cdots, a_{d}),(a_2,a_3,\cdots,a_{d},a_1),\cdots\}$$ equals a positive multiple of $kp_i$.
Then we have $$\bigcup_{i=0}^{m}\mathfrak{O}_i$$ is the set of those elements  $(a_1,a_2,\cdots,a_d)\in(\mathbb{Z}_{\geq0})^{d}$ such that the number of minimum elements in
 $$\{(a_{1},a_{2},\cdots, a_{d}),(a_2,a_3,\cdots,a_{d},a_1),\cdots\}$$ equals a positive multiple of $k$.

 Assume that $$(a_1,a_2,\cdots,a_{\frac{d}{k}})\in (\mathbb{Z}_{\geq0})^{\frac{d}{k}}$$ satisfies
$$\{a_1,a_2,\cdots,a_{\frac{d}{k}}\}=\{a_{i_1},a_{i_2},\cdots,a_{i_p}\},$$ $\#\{1\leq l\leq \frac{d}{k}|a_l=a_{i_j}\}=\frac{t_j}{k}$, $1\leq j \leq p$. If we extend it to $$(a_1,a_2,\cdots,a_{\frac{d}{k}},a_1,a_2,\cdots,a_{\frac{d}{k}},\cdots,a_1,a_2,\cdots,a_{\frac{d}{k}})=(a_1,a_2,\cdots,a_d),$$ then it is an element in $(\mathbb{Z}_{\geq0})^{d}$ such that
the number of minimum elements in
 $$\{(a_{1},a_{2},\cdots, a_{d}),(a_2,a_3,\cdots,a_{d},a_1),\cdots\}$$ equals a positive multiple of $k$, which depends on the choice of $$(a_1,a_2,\cdots,a_{\frac{d}{k}}).$$ Since $\pi_k(T,0)$ is the number of all such elements, we have $$\pi_k(T,0)=|\bigcup_{i=0}^{m}\mathfrak{O}_i|.$$

Similarly, for a fixed $\{p_{x_1},p_{x_2},\cdots, p_{x_j}\}\in P(T,k)_j$, assume that $$(a_1,a_2,\cdots,a_{\frac{d}{kp_{x_1}p_{x_2}\cdots p_{x_j}}})\in(\mathbb{Z}_{\geq0})^{\frac{d}{kp_{x_1}p_{x_2}\cdots p_{x_j}}}$$ satisfies
$$\{a_1,a_2,\cdots,a_{\frac{d}{kp_{x_1}p_{x_2}\cdots p_{x_j}}}\}=\{a_{i_1},a_{i_2},\cdots,a_{i_p}\},$$ $\#(a_{i_j})=\frac{t_j}{kp_{x_1}p_{x_2}\cdots p_{x_j}}$, $1\leq j \leq p$.  If we extend it to
$$(a_1,a_2,\cdots,a_{\frac{d}{kp_{x_1}p_{x_2}\cdots p_{x_j}}},\cdots,a_1,a_2,\cdots,a_{\frac{d}{kp_{x_1}p_{x_2}\cdots p_{x_j}}})=(a_{1},a_{2},\cdots, a_{d}),$$ then it is an element in $(\mathbb{Z}_{\geq0})^{d}$ such that
the number of minimum elements in
 $$\{(a_{1},a_{2},\cdots, a_{d}),(a_2,a_3,\cdots,a_{d},a_1),\cdots\}$$ equals a positive multiple of $kp_{x_1}p_{x_2}\cdots p_{x_j}$, which depends on the choice of $$(a_1,a_2,\cdots,a_{\frac{d}{kp_{x_1}p_{x_2}\cdots p_{x_j}}}).$$ Since the summand  of $\pi_k(T,j)$ fixed by $\{p_{x_1},p_{x_2},\cdots, p_{x_j}\}\in P(T,k)_j$ is the number of all such elements,  we have $$\pi_k(T,j)=\sum_{\{p_{x_1},p_{x_2},\cdots ,p_{x_j}\}\in P(k)_j}|\bigcap_{i=1}^{j}\mathfrak{O}_{x_i}|.$$

Since $\mathfrak{O}_0 \cap \mathfrak{O}_i=\emptyset$, $1\leq i\leq m$,
 by applying the inclusion-exclusion principle, we have
$$\pi_k(T,0)=|\mathfrak{O}_0|-\sum_{j=1}^{m}(-1)^j\pi_k(T,j).$$
Thus $$|\mathfrak{O}_0|=\sum_{j=0}^{m}(-1)^j\pi_k(T,j)$$ is the number of all
$(a_1,a_2,\cdots,a_d)\in(\mathbb{Z}_{\geq0})^{d}$ such that the number of minimum elements in
 $$\{(a_{1},a_{2},\cdots, a_{d}),(a_2,a_3,\cdots,a_{d},a_1),\cdots\}$$ equals $k$.
 Since the number of the $\langle c_{d}\rangle$-orbit of every such element equals $\frac{d}{k}$, $$\frac{k}{d}\sum_{j=0}^{m}(-1)^j\pi_k(j)$$
is the number of those invariant cycles  $\chi_i(\delta_s)=(a_1,a_2,\cdots,a_d)$ such that the number of minimum elements in
 $$\{(a_{1},a_{2},\cdots, a_{d}),(a_2,a_3,\cdots,a_{d},a_1),\cdots\}$$ equals $k$.

From Theorem $2.5$, considering all $T\in T(p)$ for all $1\leq p \leq d$, we complete our proof. \qb

From Theorem $2.7$, Theorem $3.1$ and $3.3$, since $d_1+\cdots+d_a\leq *+1$ in the range $*\leq q-1$, $H^{*}(P_n)^\mathfrak{G}$ is independent of $n$ and $q$ in the range $*\leq q-1$.

{\bf Remark:} Let $\mathfrak{G}$ be the $\mathbb{Z}/2$-extension of the direct product of two symmetric groups $\mathfrak{S}_{q} \overleftrightarrow{\times} \mathfrak{S}_q$, where the $\mathbb{Z}/2$-extension is given by $$\sigma=(1,2q)(2,2q-1)\cdots(q,q+1).$$
 Since $\mathfrak{S}_{q} \times \mathfrak{S}_q$ is a normal subgroup of
$\mathfrak{G}$, $\sigma$ acts on $H^{*}(P_n)^{\mathfrak{S}_{q} \times \mathfrak{S}_q}$ in a natural way, there exists a natural surjection
$$\phi: H^{*}(P_n)^{\mathfrak{S}_{q} \times \mathfrak{S}_q}\rightarrow H^{*}(P_n)^{\mathfrak{G}}$$
of $\mathbb{Q}\mathfrak{G}$-modules defined by $$\phi(s\otimes\zeta_{\lambda})=s\otimes\zeta_{\lambda}+\sigma(s\otimes\zeta_{\lambda}).$$
From Theorem $2.5$, Theorem $2.6$, Theorem $3.1$ and Theorem $3.3$, we already know $H^{*}(P_n)^{\mathfrak{S}_{q} \times \mathfrak{S}_q}$ clearly, hence we can calculate  $H^*(P_{n})^{\mathfrak{G}}$ by the kernel of $\phi$. The auther will give the detail in another paper.

{\noindent\bf Theorem 3.4.}
{\it Let $\mathfrak{G}=\mathfrak{S}_{n-1}\times\mathfrak{S}_1$. Then
$$H^{i}(P_n)^{\mathfrak{G}}=\left\{\begin{array}{ccc}\mathbb{Q}, & \text{if }i=0,n-1,\\
\mathbb{Q}\oplus\mathbb{Q}, & \text{if }1\leq i \leq n-2. \end{array}\right.$$}

{\noindent\bf Proof.}
Clearly we have $H^{0}(P_n)^{\mathfrak{G}}=\mathbb{Q}$.

From Lemma $2.3$,
we only need to consider the partitions
\begin{enumerate}
\item[(i)]
$\lambda=(\lambda_1,1,\cdots,1)$, $2\leq j\leq n-1$.

\item[(ii)] $\lambda=(\lambda_1,2,1,\cdots,1)$, $2\leq j\leq n-1$.

\item[(iii)] $\lambda=(2,1,\cdots,1)$.

\item[(iv)]$\lambda=(n)$.
\end{enumerate}

From Theorem $3.1$ and Theorem $3.3$, since $|\Pi(\lambda_1,1)|=1$ for $\lambda_1\geq2$, $|\Pi(2,0)|=1$,
we complete our proof. \qb

 {\noindent\bf Theorem 3.5.}
 {\it Let $\mathfrak{G}=\mathfrak{S}_{n-2}\times\mathfrak{S}_2$. Then
$$H^{i}(P_n)^{\mathfrak{G}}=\left\{\begin{array}{ccc}
\mathbb{Q}, & \text{if }i=0, &\\
\mathbb{Q}^{2i}, & \text{if }i \equiv 0,2,& 1\leq i \leq n-3 , \\
\mathbb{Q}^{2i+1}, &   \text{if }i\equiv 1,&  1\leq i\leq n-3, \\
\mathbb{Q}^{2i-1}, &  \text{if }i\equiv 3,& 1\leq i \leq n-3, \\
\mathbb{Q}^{\frac{3}{2}i}, & \text{if }i \equiv 0,2, & i=n-2,\\
\mathbb{Q}^{\frac{3}{2}i-\frac{1}{2}}, & \text{if }i \equiv 1, & i=n-2,\\
\mathbb{Q}^{\frac{3}{2}i+\frac{1}{2}}, & \text{if }i \equiv 3, & i=n-2,\\
\mathbb{Q}^{\frac{1}{2}i},  & \text{if }i \equiv 0,2, & i=n-1,\\
\mathbb{Q}^{\frac{1}{2}i+\frac{1}{2}},  & \text{if }i \equiv 1, & i=n-1,\\
\mathbb{Q}^{\frac{1}{2}i-\frac{1}{2}},  & \text{if }i \equiv 3, & i=n-1,
\end{array}\right.\text{ mod }4.$$}

{\noindent\bf Proof.}
From Lemma $2.3$, if $i\geq4$, we only need to consider the partitions
\begin{enumerate}
\item[(i)] $\lambda=(i+1,1,1,\cdots)$,  $1\leq j \leq n-4$.
\item[(ii)] $\lambda=(i,2,1,1,\cdots)$, $2\leq j \leq n-4$.
\item[(iii)] $\lambda=(\lambda_1,\lambda_2,1,1,\cdots)$, $\lambda_1\geq\lambda_2\geq 3$, $\lambda_1+\lambda_2=i+2$, $2\leq j \leq n-4$.
\item[(iv)] $\lambda=(\lambda_1,\lambda_{2},2,1,1,\cdots)$,  $\lambda_1\geq\lambda_2\geq 2$, $\lambda_1+\lambda_2=i+1$, $3 \leq j \leq n-4$.
\end{enumerate}

We already know that $|\Pi(2,2)|=|\Pi(2,1)|=|\Pi(2,0)|=1$. From Theorem $3.3$, for $\lambda_i\geq 3$, we have $|\Pi(\lambda_i,1)|=|\Pi(\lambda_i,\lambda_i-1)|=1$,
$$|\Pi(\lambda_i,2)|=
\left\{\begin{array}{cl}
\frac{1}{2}\lambda_i-1, & \text{ if } \lambda_i \equiv 0, \\
\frac{1}{2}\lambda_i-\frac{1}{2}, & \text{ if } \lambda_i \equiv 1,3, \\
\frac{1}{2}\lambda_i, & \text{ if } \lambda_i \equiv 2,
\end{array}\right. \text{ mod } 4.
$$

Now we apply Theorem $3.1$ for $1\leq i\leq n-1$.
\begin{enumerate}
\item When $i=1$, the dimension of  $H^{1}(P_n)^{\mathfrak{G}}$  equals $$|\Pi(2,2)|+|\Pi(2,1)|+|\Pi(2,0)|=3.$$

\item When $i=2$, we only need to consider $\lambda=(3,1,1,\cdots)$ and $\lambda=(2,2,1,\cdots)$.

When $\lambda=(3,1,1,\cdots)$,  the dimension of the contribution of this case to $H^{2}(P_n)^{\mathfrak{G}}$  equals
$$|\Pi(3,1)|+|\Pi(3,2)|=2.$$

  When $\lambda=(2,2,1,\cdots)$,  the dimension  of the contribution of this case to $H^{2}(P_n)^{\mathfrak{G}}$  equals
  $$|\Pi(2,2)||\Pi(2,0)|+|\Pi(2,1)||\Pi(2,0)|=2.$$

  Hence $H^{2}(P_n)^{\mathfrak{G}}=\mathbb{Q}^4$.

\item When $i=3$, we only need to consider $\lambda=(4,1,1,\cdots)$, $\lambda=(3,2,1,\cdots)$ and $\lambda=(2,2,2,1,1,\cdots)$.

When $\lambda=(4,1,1,\cdots)$,  the dimension  of the contribution of this case to $H^{3}(P_n)^{\mathfrak{G}}$  equals
$$|\Pi(4,1)|+|\Pi(4,2)|=2.$$

  When $\lambda=(3,2,1,\cdots)$, the dimension  of the contribution of this case to $H^{3}(P_n)^{\mathfrak{G}}$  equals
  $$|\Pi(3,2)||\Pi(2,0)|+|\Pi(3,1)||\Pi(2,0)|+|\Pi(3,1)||\Pi(2,1)|=3.$$

  When $\lambda=(2,2,2,1,\cdots)$,
  the dimension  of the contribution of this case to $H^{3}(P_n)^{\mathfrak{G}}$  equals $0$.

  Hence $H^{3}(P_n)^{\mathfrak{G}}=\mathbb{Q}^5$.

\item When $4\leq i \leq n-3$,
  the dimension  of the contribution of the case (i) to $H^{i}(P_n)^{\mathfrak{G}}$  equals
   $$|\Pi(i+1,1)|+|\Pi(i+1,2)|.$$

The dimension of the contribution of the case (ii) to $H^{i}(P_n)^{\mathfrak{G}}$  equals
$$|\Pi(i,2)||\Pi(2,0)|+|\Pi(i,1)||\Pi(2,1)|+|\Pi(i,1)||\Pi(2,0)|.$$

For the case (iii), there are $\lfloor\frac{i-2}{2}\rfloor$ partitions for this case.
We remark that when $\lambda=(2k,2k,1,1,\cdots)$, $$\left(\begin{array}{c}|\Pi(2k,1)|\\ 2\end{array}\right)=0.$$
Hence
the dimension of the contribution of these partitions to $H^{i}(P_n)^{\mathfrak{G}}$  equals
\begin{align*}
\begin{array}{cl}
\frac{i-4}{2}, &\text{ if } i \equiv0,\\
\frac{i-3}{2}, & \text{ if } i \equiv1,
 \end{array} \text{ mod } 2.
 \end{align*}

For the case (iv), there are $\lfloor\frac{i-3}{2}\rfloor$ partitions for this case.
We remark that when $\lambda=(2k,2k,2,1,\cdots)$, $$\left(\begin{array}{c}|\Pi(2k,1)|\\ 2\end{array}\right)=0.$$
Hence
the dimension of the contribution of these partitions to $H^{i}(P_n)^{\mathfrak{G}}$ equals
$$
\begin{array}{cl}
\frac{i-5}{2}, & \text{ if } i-1 \equiv 0,\\
\frac{i-4}{2}, &\text{ if } i-1 \equiv1,
\end{array} \text{ mod } 2.
$$

\item When $i=n-2$, we only need to consider the case (i), (ii) and (iii).

 The dimension of the contribution of the case (i) and (iii) to $H^{i}(P_n)^{\mathfrak{G}}$  are the same with the case $(4)$, (i) and (iii).

The dimension of the contribution of the case (ii) to $H^{i}(P_n)^{\mathfrak{G}}$  equals
$$|\Pi(i,2)||\Pi(2,0)|+|\Pi(i,1)||\Pi(2,1)|.$$

\item When $i=n-1$,
we only need to consider the case (i).
 The dimension of the contribution of the case (i) to $H^{i}(P_n)^{\mathfrak{G}}$  equals
   $$|\Pi(i+1,2)|.$$
\end{enumerate}

Hence by calculating $i$ mod $4$, we complete our proof. \qb

{\noindent\bf Theorem 3.6.}
{\it Let $\mathfrak{G}=\mathfrak{S}_{n-3}\times\mathfrak{S}_3$. Then
$$H^{i}(P_n)^{\mathfrak{G}}=\left\{\begin{array}{ccc}
\mathbb{Q}^{i^2-\frac{1}{3}i+1}, & \text{if }i \equiv 0,3, & 0\leq i\leq n-4,\\
\mathbb{Q}^{i^2-\frac{1}{3}i+\frac{7}{3}}, & \text{if }i \equiv 1,10, & 0\leq i\leq n-4,\\
\mathbb{Q}^{i^2-\frac{1}{3}i+\frac{5}{3}}, & \text{if }i \equiv 2,5, & 0\leq i\leq n-4,\\
\mathbb{Q}^{i^2-\frac{1}{3}i+\frac{4}{3}},  & \text{if }i \equiv 4,7, & 0\leq i\leq n-4,\\
\mathbb{Q}^{i^2-\frac{1}{3}i+2,}  & \text{if }i \equiv 6,9, &0\leq i\leq n-4\\
\mathbb{Q}^{i^2-\frac{1}{3}i+\frac{2}{3}},  & \text{if }i \equiv 8,11, & 0\leq i\leq n-4,\\
\mathbb{Q}^{\frac{11}{12}i^2-\frac{5}{12}i+1},  & \text{if }i\equiv 0,3,4,7, &i=n-3,\\
\mathbb{Q}^{\frac{11}{12}i^2-\frac{5}{12}i+\frac{3}{2}},  & \text{if }i\equiv 1,6,9,10, &i=n-3,\\
\mathbb{Q}^{\frac{11}{12}i^2-\frac{5}{12}i+\frac{7}{6}},  & \text{if }i\equiv 2,5,&i=n-3,\\
\mathbb{Q}^{\frac{11}{12}i^2-\frac{5}{12}i+\frac{2}{3}},  & \text{if }i\equiv 8,11,&i=n-3,\\
\mathbb{Q}^{\frac{7}{12}i^2-\frac{5}{12}i},  & \text{if }i\equiv 0,3,8,11,&i=n-2,\\
\mathbb{Q}^{\frac{7}{12}i^2-\frac{5}{12}i+\frac{5}{6}},  & \text{if }i\equiv 1,10,&i=n-2,\\
\mathbb{Q}^{\frac{7}{12}i^2-\frac{5}{12}i+\frac{1}{2}},  & \text{if }i\equiv 2,5,6,9,&i=n-2,\\
\mathbb{Q}^{\frac{7}{12}i^2-\frac{5}{12}i+\frac{1}{3}},  & \text{if }i\equiv 4,7,&i=n-2,\\
\mathbb{Q}^{\frac{1}{6}i^2-\frac{1}{6}i},  & \text{if }i \equiv 0,1,3,4,6,7,9,10,& i=n-1,\\
\mathbb{Q}^{\frac{1}{6}i^2-\frac{1}{6}i-\frac{1}{3}},  & \text{if }i \equiv 2,5,8,11,& i=n-1,\\
\end{array}\right.\text{ mod }12.$$}

{\noindent\bf Proof.}
From Lemma $2.3$, if $i\geq5$, we only need to consider the partitions
\begin{enumerate}
\item[(i)] $\lambda=(i+1,1,1,\cdots)$,  $1\leq j \leq n-5$.
\item[(ii)] $\lambda=(i,2,1,1,\cdots)$, $2\leq j \leq n-5$.
\item[(iii)] $\lambda=(\lambda_1,\lambda_2,1,1,\cdots)$, $\lambda_1\geq\lambda_2\geq 3$, $\lambda_1+\lambda_2=i+2$, $2\leq j \leq n-5$.
\item[(iv)] $\lambda=(\lambda_1,\lambda_{2},2,1,1,\cdots)$,  $\lambda_1\geq\lambda_2\geq 2$, $\lambda_1+\lambda_2=i+1$, $3 \leq j \leq n-5$.
\item[(v)] $\lambda=(\lambda_1,\lambda_{2},\lambda_{3},1,1,\cdots)$,  $\lambda_1\geq\lambda_2\geq \lambda_{3}\geq3$, $\lambda_1+\lambda_2+\lambda_{3}=i+3$, $3 \leq j \leq n-5$.
\item[(vi)] $\lambda=(\lambda_1,\lambda_{2},\lambda_{3},2,1,\cdots)$,  $\lambda_1\geq\lambda_2\geq \lambda_{3}\geq2$, $\lambda_1+\lambda_2+\lambda_{3}=i+2$, $4 \leq j \leq n-5$.
\end{enumerate}

 We already know that $|\Pi(2,2)|=|\Pi(2,1)|=|\Pi(2,0)|=1$.
From Theorem $3.3$ and the proof of Theorem $3.5$, for $\lambda_i\geq 3$,
$$|\Pi(\lambda_i,2)|=
\left\{\begin{array}{cl}
\frac{1}{2}\lambda_i-1, & \text{ if } \lambda_i \equiv 0, \\
\frac{1}{2}\lambda_i-\frac{1}{2}, & \text{ if } \lambda_i \equiv 1,3, \\
\frac{1}{2}\lambda_i, & \text{ if } \lambda_i \equiv 2,
\end{array}\right. \text{ mod } 4.
$$
$$|\Pi(\lambda_i,3)|=
\left\{\begin{array}{cl}
\frac{1}{6}\lambda_{i}^2-\frac{1}{2}\lambda_i , & \text{if } \lambda_i \equiv 0, \\
\frac{1}{6}\lambda_{i}^2-\frac{1}{2}\lambda_i+\frac{1}{3} , & \text{if } \lambda_i \equiv 1,2,
\end{array}\right. \text{ mod } 3.
$$

Now we apply Theorem $3.1$ for $1\leq i \leq n-1$.

\begin{enumerate}
\item When $i=1$, from Theorem $2.7$ and Theorem $3.4$, we have  $H^{1}(P_n)^{\mathfrak{G}}=\mathbb{Q}^3$.

\item When $i=2$, we only need to consider the partitions $\lambda=(3,1,1,\cdots)$ and $\lambda=(2,2,1,\cdots)$.

When $\lambda=(3,1,1,\cdots)$, the dimension of the contribution of this case to $H^{2}(P_n)^{\mathfrak{G}}$  equals
$$|\Pi(3,1)|+|\Pi(3,2)|=2.$$

When $\lambda=(2,2,1,\cdots)$,the dimension of the contribution of this case to $H^{2}(P_n)^{\mathfrak{G}}$  equals
  $$|\Pi(2,2)||\Pi(2,1)|+|\Pi(2,2)||\Pi(2,0)|+|\Pi(2,1)||\Pi(2,0)|=3.$$

 Hence we have $H^{2}(P_n)^{\mathfrak{G}}=\mathbb{Q}^5$.

\item When $i=3$, we only need to consider $\lambda=(4,1,1,\cdots)$, $\lambda=(3,2,1,\cdots)$ and $\lambda=(2,2,2,1,1,\cdots)$.

When $\lambda=(4,1,1,\cdots)$,  the dimension of the contribution of this case to $H^{3}(P_n)^{\mathfrak{G}}$  equals
$$|\Pi(4,1)|+|\Pi(4,2)|+|\Pi(4.3)|=3.$$

  When $\lambda=(3,2,1,\cdots)$, the dimension of the contribution of this case to $H^{3}(P_n)^{\mathfrak{G}}$  equals
  $$\left(\begin{array}{c}|\Pi(3,2)||\Pi(2,0)|+|\Pi(3,2)||\Pi(2,1)|\\+|\Pi(3,1)||\Pi(2,0)|+|\Pi(3,1)||\Pi(2,1)|+|\Pi(3,1)||\Pi(2,2)|
  \end{array}\right)=5.$$

  When $\lambda=(2,2,2,1,\cdots)$,
  the dimension of the contribution of this case to $H^{3}(P_n)^{\mathfrak{G}}$ equals
  $$|\Pi(2,0))||\Pi(2,1)||\Pi(2,2)|=1.$$

 Hence we have $H^{3}(P_n)^{\mathfrak{G}}=\mathbb{Q}^9$.

 \item When $i=4$, we need to consider the partitions $\lambda=(5,1,1,\cdots)$, $\lambda=(4,2,1,\cdots)$, $\lambda=(3,3,1,\cdots)$,  $\lambda=(3,2,2,1,\cdots)$ and $\lambda=(2,2,2,2,1,\cdots)$.

When $\lambda=(5,1,1,\cdots)$, the dimension of the contribution of this case to $H^{4}(P_n)^{\mathfrak{G}}$  equals
$$|\Pi(5,1)|+|\Pi(5,2)|+|\Pi(5.3)|=5.$$

When $\lambda=(4,2,1,\cdots)$, the dimension of the contribution of this case to $H^{4}(P_n)^{\mathfrak{G}}$  equals
$$\left(\begin{array}{r}
|\Pi(4,3)||\Pi(2,0)|+|\Pi(4,2)||\Pi(2,0)|+|\Pi(4,2)||\Pi(2,1)|\\+|\Pi(4,1)||\Pi(2,0)|+|\Pi(4,1)||\Pi(2,1)|+|\Pi(4,1)||\Pi(2,2)|
  \end{array}\right)
  =6.
  $$

When $\lambda=(3,3,1,\cdots)$, the dimension of the contribution of this case to $H^{4}(P_n)^{\mathfrak{G}}$  equals
$$
|\Pi(3,2)||\Pi(3,1)|+\left(\begin{array}{c}|\Pi(3,1)|\\1\end{array}\right) =2.
  $$

When $\lambda=(3,2,2,1,\cdots)$, the dimension of the contribution of this case to $H^{4}(P_n)^{\mathfrak{G}}$  equals
$$\left(\begin{array}{r}
|\Pi(3,2)||\Pi(2,1)||\Pi(2,0)|\\+|\Pi(3,1)||\Pi(2,1)||\Pi(2,0)|\\+|\Pi(3,1)||\Pi(2,2)||\Pi(2,0)|  \end{array}\right)=3.
  $$

When $\lambda=(2,2,2,2,1,\cdots)$, the dimension of the contribution of this case to $H^{4}(P_n)^{\mathfrak{G}}$  equals $0$.

Hence we have $H^{4}(P_n)^{\mathfrak{G}}=\mathbb{Q}^{16}$.

\item When $5\leq i\leq n-4$, the dimension of the contribution of the case (i) to $H^{i}(P_n)^{\mathfrak{G}}$  equals
   $$|\Pi(i+1,1)|+|\Pi(i+1,2)|+|\Pi(i+1,3)|.$$

    The dimension of the contribution of the case (ii) to $H^{i}(P_n)^{\mathfrak{G}}$  equals
   $$\left(\begin{array}{r}|\Pi(i,3)||\Pi(2,0)|+|\Pi(i,2)||\Pi(2,0)|+|\Pi(i,2)||\Pi(2,1)|\\
   +|\Pi(i,1)||\Pi(2,0)|+|\Pi(i,1)||\Pi(2,1)|+|\Pi(i,1)||\Pi(2,2)|\end{array}\right).$$

For the case (iii), there are $\lfloor\frac{i-2}{2}\rfloor$ partitions for this case.
We remark that when $\lambda=(2k,2k,1,1,\cdots)$, $$\left(\begin{array}{c}|\Pi(2k,1)|\\ 2\end{array}\right)=0.$$
Hence
the dimension of the contribution of these partitions to $H^{i}(P_n)^{\mathfrak{G}}$ equals
$$
\sum_{3\leq j\leq i-1}(|\Pi(j,2)||\Pi(i+2-j,1)|)
+\left\{\begin{array}{cl}
\frac{i-4}{2}, &\text{ if } i \equiv0,\\
\frac{i-3}{2}, & \text{ if } i \equiv1,
 \end{array}\right. \text{ mod } 2.
$$

For the case (iv), there are $\lfloor\frac{i-3}{2}\rfloor$ partitions for this case.
We remark that when $\lambda=(2k,2k,2,1,\cdots)$, $$\left(\begin{array}{c}|\Pi(2k,1)|\\ 2\end{array}\right)=0.$$
Hence
the dimension of the contribution of these partitions to $H^{i}(P_n)^{\mathfrak{G}}$  equals
$$
\sum_{3\leq j\leq i-2}(|\Pi(j,2)||\Pi(i+1-j,1)|)+\left\{\begin{array}{cl}
i-5, & \text{ if } i-1 \equiv 0,\\
i-4, &\text{ if } i-1 \equiv1,
\end{array}\right. \text{ mod } 2.
$$

For the case (v), there are
$$
\begin{array}{cl}
\frac{1}{12}i^2-\frac{1}{2}i+1, & \text{ if } i \equiv 0,\\
\frac{1}{12}i^2-\frac{1}{2}i+\frac{5}{12}, &\text{ if } i \equiv1,5,\\
\frac{1}{12}i^2-\frac{1}{2}i+\frac{2}{3},  &\text{ if } i \equiv2,4,\\
\frac{1}{12}i^2-\frac{1}{2}i+\frac{3}{4},  &\text{ if } i \equiv3,
\end{array} \text{ mod } 6,
$$
 partitions for this case.
We remark that when $$\lambda=(\lambda_1,2k,2k,1,\cdots)$$ or $$\lambda=(2k,2k,\lambda_3,1,\cdots),$$  $\lambda_1\geq 2k\geq \lambda_3\geq2$,  $$\left(\begin{array}{c}|\Pi(2k,1)|\\ 2\end{array}\right)=\left(\begin{array}{c}|\Pi(2k,1)|\\ 3\end{array}\right)=0.$$ Totally there are $\lfloor\frac{i-4}{4}\rfloor$ such partitions.
Hence
the dimension of the contribution of these partitions to $H^{i}(P_n)^{\mathfrak{G}}$  equals
$$
\begin{array}{cl}
\frac{1}{12}i^2-\frac{3}{4}i+2, & \text{ if } i \equiv 0,9 ,\\
\frac{1}{12}i^2-\frac{3}{4}i+\frac{5}{3}, &\text{ if } i \equiv1,4,5,8,\\
\frac{1}{12}i^2-\frac{3}{4}i+\frac{13}{6},  &\text{ if } i \equiv2,7,10,11,\\
\frac{1}{12}i^2-\frac{3}{4}i+\frac{5}{2},  &\text{ if } i \equiv3,6,
\end{array} \text{ mod } 12.
$$

For the case (vi), there are
$$
\begin{array}{cl}
\frac{1}{12}i^2-\frac{1}{6}i, & \text{ if } i \equiv 0,2,\\
\frac{1}{12}i^2-\frac{1}{6}i+\frac{1}{12}, &\text{ if } i \equiv1,\\
\frac{1}{12}i^2-\frac{1}{6}i-\frac{1}{4},  &\text{ if } i \equiv3,5,\\
\frac{1}{12}i^2-\frac{1}{6}i+\frac{1}{3},  &\text{ if } i \equiv4,
\end{array} \text{ mod } 6,
$$
 partitions for this case.
We remark that when $$\lambda=(\lambda_1,2k,2k,2,1,\cdots)$$ or $$\lambda=(2k,2k,\lambda_3,2,1,\cdots),$$  $\lambda_1\geq 2k\geq \lambda_3\geq2$,  $$\left(\begin{array}{c}|\Pi(2k,1)|\\ 2\end{array}\right)=\left(\begin{array}{c}|\Pi(2k,1)|\\ 3\end{array}\right)=0.$$ Totally there are $\lfloor\frac{i}{4}\rfloor$ such partitions.
Hence
the dimension of the contribution of these partitions to $H^{i}(P_n)^{\mathfrak{G}}$  equals
$$
\begin{array}{cl}
\frac{1}{12}i^2-\frac{5}{12}i, & \text{ if } i \equiv 0,5,8,9,\\
\frac{1}{12}i^2-\frac{5}{12}i+\frac{1}{3}, &\text{ if } i \equiv1,4,\\
\frac{1}{12}i^2-\frac{5}{12}i+\frac{1}{2}, &\text{ if } i \equiv2,3,6,11,\\
\frac{1}{12}i^2-\frac{5}{12}i+\frac{5}{6}, & \text{ if } i \equiv 7,10,
\end{array} \text{ mod } 12.
$$

\item When $i=n-3$,
we only need to consider the case (i), (ii), (iii), (iv) and (v).

The dimension of the contribution of the case (i), (iii) and (v) to $H^{i}(P_n)^{\mathfrak{G}}$  are the same with the case $(4)$, (i), (iii) and (v).

The dimension of the contribution of the case (ii) to $H^{i}(P_n)^{\mathfrak{G}}$ equals
 $$\left(\begin{array}{c}|\Pi(i,3)||\Pi(2,0)|+|\Pi(i,2)||\Pi(2,0)|+|\Pi(i,2)||\Pi(2,1)|\\
   +|\Pi(i,1)||\Pi(2,1)|+|\Pi(i,1)||\Pi(2,2)|\end{array}\right).$$

The dimension of the contribution of the case (iv) to $H^{i}(P_n)^{\mathfrak{G}}$  equals
$$\sum_{3\leq j\leq i-2}(|\Pi(j,2)||\Pi(i+1-j,1)|)+\left\{\begin{array}{cl}
\frac{i-5}{2}, & \text{ if } i-1 \equiv 0,\\
\frac{i-4}{2}, &\text{ if } i-1 \equiv1,
\end{array}\right. \text{ mod } 2.
$$

\item When $i=n-2$,
we only need to consider the case (i), (ii) and (iii).

The dimension of  the contribution of the case (iii) to $H^{i}(P_n)^{\mathfrak{G}}$  are the same with the case $(4)$, (iii).

The dimension of the contribution of the case (i) to $H^{i}(P_n)^{\mathfrak{G}}$ equals
$$|\Pi(i+1,2)|+|\Pi(i+1,3)|.$$

The dimension of the contribution of the case (ii) to $H^{i}(P_n)^{\mathfrak{G}}$ equals
$$|\Pi(i,3)||\Pi(2,0)|+|\Pi(i,2)||\Pi(2,1)|+|\Pi(i,1)||\Pi(2,2)|.$$

\item When $i=n-1$, we only need to consider the case (i). The dimension of the contribution of the case (i) to $H^{i}(P_n)^{\mathfrak{G}}$ equals
   $$|\Pi(i+1,3)|.$$
\end{enumerate}

Hence by calculating $i$ mod $12$, we complete our proof. \qb

\end{document}